\documentclass[a4paper,twoside,12pt]{amsart}

\usepackage{graphicx,amsmath,amssymb,hyperref}

\newenvironment{proof2}{\begin{proof}}{\qed\end{proof}}
\newenvironment{acknowledgements}{\paragraph {\bf Ackowledgements}}{}

\def \Res{\operatorname{Res}}     
\def \sr{\operatorname{sr}}       
\def \P{\mathbb{P}}               
\def \RP{\mathbb{RP}}             
\def \C{{\mathbb{C}}}             
\def \R{{\mathbb{R}}}             
\def \Disc{\operatorname{Disc}}   
\def \Diag{\operatorname{Diag}}   
\def \id{\operatorname{id}}

\def \PI{\textsf{I}}
\def \PIa{\textsf{Ia}}
\def \PIb{\textsf{Ib}}
\def \PII{\textsf{II}}
\def \PIIa{\textsf{IIa}}
\def \PIII{\textsf{III}}
\def \PIIIa{\textsf{IIIa}}
\def \PIV{\textsf{IV}}
\def \PV{\textsf{V}}

\def \CIN{\textsf{IN}}
\def \CIS{\textsf{IS}}
\def \CIaN{\textsf{IaN}}
\def \CIaS{\textsf{IaS}}
\def \CIbN{\textsf{IbN}}
\def \CIIN{\textsf{IIN}}
\def \CIIS{\textsf{IIS}}
\def \CIIaS{\textsf{IIaS}}
\def \CIIaN{\textsf{IIaN}}
\def \CIIIN{\textsf{IIIN}}
\def \CIIIS{\textsf{IIIS}}
\def \CIIIaN{\textsf{IIIaN}}
\def \CIVN{\textsf{IVN}}
\def \CVN{\textsf{VN}}

\def \CN{\textsf{N}}
\def \CS{\textsf{S}}

\newtheorem{proposition}{\textsc{Proposition}}
\newtheorem{definition}{\textsc{Definition}}
\newtheorem{theorem}{\textsc{Theorem}}
\newtheorem{corollary}{\textsc{Corollary}}
\newtheorem{lemma}{\textsc{Lemma}}

\newenvironment{remark}{\paragraph*{\bf Remark~:}}{\hfill $\boxdot$ \smallskip}

\newcommand{\aaeccarray}[3]{
\begin{array}{#1} 
\hline
#2\\
\hline
#3
\hline
\end{array}
}

\title[The configurations of two real projective conics]{Equations, inequations and inequalities characterizing the
  configurations of two real projective conics}

\date{November 2005}

\author{Emmanuel Briand}
\address{
Emmanuel Briand\\
Universidad de Cantabria\\
Dpto. Matem\'aticas, estad\'{\i}stica y computaci\'on\\
Avda. Los Castros S/N\\
39005 Santander\\
Spain.
}
\email{ebriand@us.es}
\urladdr{http://emmanuel.jean.briand.free.fr/}

\begin{document}
\bibliographystyle{plain}

\begin{abstract}
Couples of 
proper, non-empty real
projective conics can be classified \emph{modulo} rigid isotopy and ambient isotopy.
We characterize the classes 
by equations, inequations and 
inequalities in the coefficients of the quadratic forms defining the
conics.
The results are well--adapted to the study of the relative position of two
conics defined by equations depending on parameters.
\keywords{arrangements of conics, rigid isotopy, relative
position of two conics, classical invariant theory.}
\end{abstract}

\maketitle

{\bf MSC2000:} 13A50 (invariant theory), 13J30 (real algebra).

\smallskip

{\emph{This is the  preliminary version of the paper  published in: Applicable Algebra in Engineering, Communication and Computing vol. 18 (1-2), pp. 21-52 (2007). This is not the version of record, which is available at} 
\url{https://dx.doi.org/10.1007/s00200-006-0023-8}.}



\section{Introduction}\label{section:introduction}


Couples of proper real projective conics, admitting
real points, can be classified \emph{modulo} ambient isotopy.
The goal of this paper is to provide equations, inequations and
inequalities characterizing each class.

This is
particularly well-suited for the following problem: given two conics
whose equations depend on parameters, for which values of the
parameters are these conics in a given ambient isotopy class ? 

Such problems are of interest in geometric modeling. They are considered for instance in the articles
\cite{E:GV:dR,Wang:Krasauskas} (and \cite{Wang:Wang:Kim} for the
similar problem for ellipsoids). We consider this paper
as the systematization of their main ideas.
Specially, in
\cite{E:GV:dR}, an algorithm was proposed to determine the configuration of a
pair of ellipses, by means of calculations of Sturm-Habicht
sequences. Our approach is different: when there, computations were
performed for each particular case, we perform the computations once
for all
in the
most general case. The formulas obtained behave well under specialization.

\begin{figure}
\centering
      \[
      \begin{array}{ccc}
	     \multicolumn{3}{c}{ 
                \begin{array}{lr}
		  \includegraphics[height=1.3cm]{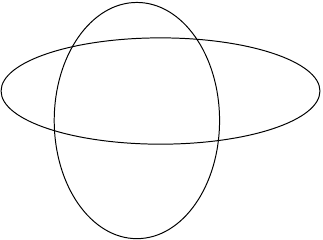} &
		  \includegraphics[height=1.3cm]{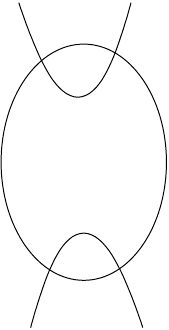}  \\ 
	     {{\CIN}} & {\CIS} 
       		\end{array}
	     }\\
	     &&\\
	       \includegraphics[height=1.3cm]{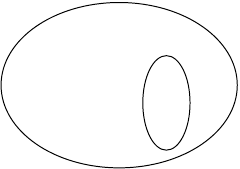} &
	       \includegraphics[height=1.3cm]{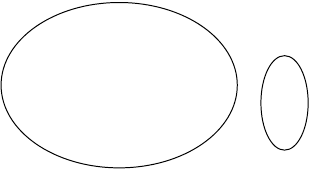} &
	       \includegraphics[height=1.3cm]{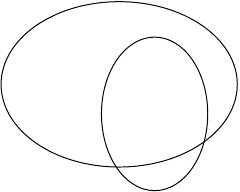}\\
	       {\CIaN} (*) & {\CIaS} & {\CIbN}   
     \end{array}
\]
\caption{The rigid isotopy classes for generic pairs of conics.}\label{fig:generic}
\end{figure}

\begin{figure}
\centering
\[
\begin{array}{cccc}
         \includegraphics[height=1.3cm]{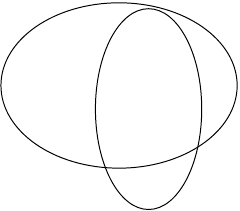} &
         \includegraphics[height=1.3cm]{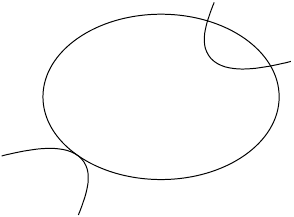} &
         \includegraphics[height=1.3cm]{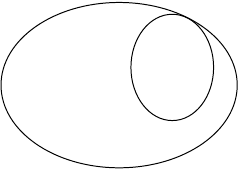} &
         \includegraphics[height=1.3cm]{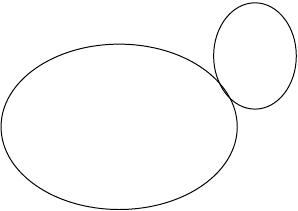} \\
	 {\CIIN} (*) & {\CIIS} & {\CIIaN} (*) & {\CIIaS} \\
         &&&\\
	 \multicolumn{4}{c}{
	   \begin{array}{ccc}
             \includegraphics[height=1.3cm]{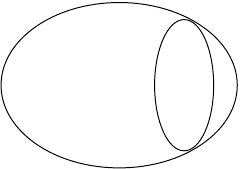} &
             \includegraphics[height=1.3cm]{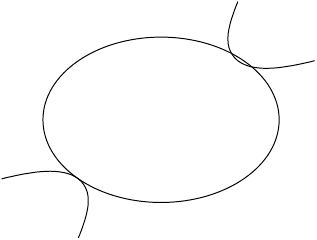} &
             \includegraphics[height=1.3cm]{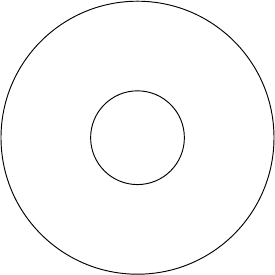} \\
	      {\CIIIN} (*) & {\CIIIS} & {\CIIIaN} (*) 
	   \end{array}  
	 }\\
         &&&\\
	 \multicolumn{4}{c}{
	   \begin{array}{cc}
             \includegraphics[angle=90,height=1.3cm]{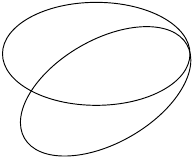} &
             \includegraphics[angle=90,height=1.3cm]{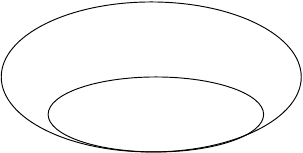} \\
	      {\CIVN} & {\CVN} (*) 
	   \end{array}  
	 }\\
\end{array}
\]
\caption{The rigid isotopy classes for non-generic pairs of proper conics.}\label{fig:singular}
\end{figure}

Instead of working with ambient isotopy, we consider another
equivalence relation, \emph{rigid isotopy}\footnote{We follow the
  terminology used in real algebraic geometry in similar situations.}, corresponding to real
deformation of the equations of the conics that doesn't change the
nature of the (complex) singularities (definition
\ref{definition:rigid_isotopy}, following the ideas of \cite{Gudkov}). Figures \ref{fig:generic} and \ref{fig:singular} provide a drawing
for a representative of each class.
Rigid isotopy happens to be an
equivalence relation just slightly finer than ambient isotopy. So we
get the classification under ambient isotopy directly from the one
under rigid isotopy.

The classification of pairs of real projective conics under rigid
isotopy was first obtained
by Gudkov and Polotovskiy
\cite{Gudkov:Polotovskiy:1,Gudkov:Polotovskiy:2,Gudkov:Polotovskiy:3}
in their work on quartic real projective
curves. Nevertheless, we start (section 2) with re--establishing this
classification. We emphasize the following key ingredient: that 
any rigid
isotopy decomposes into a path in one orbit of the space of
pencils of conics under projective transformations, and a
rigid isotopy stabilizing some pencil of conic (Lemma \ref{lemma:basic_ri}). 
As a consequence, each rigid isotopy class is determined by an orbit
of pencils of conics and the position of the two conics with respect to the degenerate
conics in the pencil they generate. 
This has two direct applications. 

First, because there are finitely many
orbits of pencils of conics under projective transformations, we get easily a finite
set of couples of conics meeting at least once each rigid isotopy
class. 

Second, it indicates clearly how to derive the equations, inequations,
inequalities characterizing the classes, which is done in section 3. The determination of
the 
position of the conics with respect to the degenerate conics in a pencil
essentially reduces to 
problems of location of roots of univariate polynomials. They can be
treated using standard tools from real algebra, namely Descartes' law
of signs and subresultant sequences. This contributes both to 
the characterization of the orbits of pencils of conics, and
the characterization of the rigid isotopy class associated to each
orbit of pencils of conics.
Classical invariant theory is also used for the first task.

Last, section 4 provides some examples of computations using the
previous results.

\subsection*{Generalities and notations}

The real projective space of dimension $k$ will be denoted with
$\R\P^k$; in particular, $\RP^2$ denotes the projective plane. 
The space of real ternary quadratic forms will be denoted with $S^2
{\R^3}^*$.
We will consider $\P(S^2 {\R^3}^*)$, the associated projective space
(see \cite{CLO} for the definitions of the notions of projective
geometry needed here).
The term \emph{conic} will be used with two meanings:
\begin{itemize}
\item an algebraic meaning: an element of  $\P(S^2
  {\R^3}^*)$. The algebraic conic associated to the quadratic form $f$
  will be denoted with $[f]$.
\item a geometric meaning: the zero locus, in $\R\P^2$, of a
  non-zero quadratic form $f$. It will be denoted with $[f=0]$.
\end{itemize}
A (geometric or algebraic) conic is said \emph{proper} if
it comes from a non-degenerate quadratic form; \emph{degenerate} if it comes
from a degenerate quadratic form. Note that, with this definition, the
empty set is a proper (geometric) conic. Algebraic and geometric
proper non-empty conics are in bijection, and can be identified. 

We define the \emph{discriminant of the quadratic form $f$} to be
\[
\Disc(f)
=
\det(\operatorname{Matrix}(f)).
\]

Any proper non-empty conic cuts out the real projective
plane into two connected components. 
They are topologically
non-equivalent: one is homeomorphic to a M\"obius strip, the other to an
open disk. The former is the \emph{outside} of the conic, the latter
is its \emph{inside}.
Let $f_0=x^2 + y^2 - z^2$.
The inside of $[f_0=0]$ is the solution set of the inequation
$f_0 < 0$, or, equivalently, the set of
points where $f_0$ has the sign of $\Disc(f_0)$. These signs change
together under linear transformations.
Now any proper, non-empty conic is obtained from $[f_0=0]$ by means of a
transformation of $PGL(3,\R)$.
Thus \emph{the inside of $[f=0]$ is the set of points where $f$ takes the sign
of $\Disc(f)$}.

The \emph{tangential quadratic form} associated to the quadratic form
$f$ on $\R^3$
is the quadratic form $\tilde{f}$ on ${\R^3}^*$ whose matrix is the
matrix of the cofactors of the matrix of $f$. The $\emph{tangential
  conic}$ associated to $[f]$ (resp. $[f=0]$) is $[\tilde{f}]$
(resp. $[\tilde{f}=0]$).

A \emph{pencil of quadratic forms} is a plane (through the origin) in ${S^2 \R^3}^*$; the
associated \emph{(projective) pencil of conics} is the corresponding
line in $\P({S^2 \R^3}^*)$. It is said to be \emph{non-degenerate} if
it contains proper conics\footnote{Contrary, for instance, to the
  pencil of the zero loci of the $f(x,y,z)=\lambda x y + \mu x z$.}.
The common points of all conics of a given pencil are called the
\emph{base points} of the pencil. They are also the common points of
any two distinct conics of the pencil. 
A non-degenerate pencil of conics has always four common points in the complex
projective space, when counted with multiplicities. 

Note that we will distinguish between (ordered) \emph{couples} of
conics ($(C_1,C_2)$ distinct from $(C_2,C_1)$, except when $C_1=C_2$) and
(unordered) \emph{pairs} of conics ($\{C_1,C_2\}=\{C_2,C_1\}$).

The \emph{characteristic form} of the couple $(f,g)$ of real ternary
quadratic forms is the binary cubic in $(t,u)$:
\[
\Phi(f,g;t,u)
:=
\Disc(t f + u g).
\]
Its coefficients will be denoted as follows:
\[
\Phi(f,g;t,u)
=
\Phi_{30} t^3 
+\Phi_{21} t^2 u
+\Phi_{12} t u^2 
+\Phi_{30} u^3. 
\]
We will also consider the de-homogenized polynomial obtained from
$\Phi$ by setting $u=1$. It will be denoted with $\phi(f,g;t)$, or $\phi(t)$
when there is no ambiguity about $f,g$. So:
\[
\phi(t):=\Disc(t f + g)
\]
Note that $\Phi_{30}=\Disc(f)$ and $\Phi_{03}=\Disc(g)$.

An \emph{isotopy} of a manifold $M$ is a continuous mapping $\theta: I\times M
\rightarrow M$, where $I$ is an interval containing $0$, such that for each $t\in I$,
the mapping $x \mapsto \theta(t,x)$ is an homeomorphism of $M$ onto
itself, and $x\mapsto \theta(0,x)$ is the identity of $M$.

Two subsets $N_1,N_2$ of $M$ are \emph{ambient isotopic} if there is an
isotopy of $M$ such that, at some instant $t\in I$, $\theta(t,N_1)=N_2$.

This definition is immediately generalized to couples of subsets:
$(N_1,N'_1)$ and $(N_2,N'_2)$ are ambient isotopic if there is an
isotopy of $M$ such that, at some instant $t$, $\theta(t,N_1)=N_2$ and $\theta(t,N'_1)=N'_2$.

\section{Classification}

\subsection{Rigid isotopy}

To classify the couples of conics up to ambient isotopy, we introduce a
slightly finer equivalence relation, \emph{rigid isotopy},
corresponding to a continuous path in the space of couples of
distinct proper conics, that doesn't change the nature of the complex
singularities of the union 
of the conics. Before stating formally the definition (definition
\ref{definition:rigid_isotopy} below), we clarify this point.
The complex singularities of
the union of the conics correspond to the (real and imaginary)
intersections of the conics. For a given multiplicity, there is
only one analytic type of intersection point of two conics. Thus the
nature of the singularities for the union of two distinct proper conics is
determined by the numbers of real and imaginary intersections of each
multiplicity. This is narrowly connected to the projective classification of
pencils of conics, that can be found in \cite{Degtyarev,Levy}. The
connection is the following theorem.
\begin{theorem}\label{pencils:base_points}\emph{(\cite{Degtyarev,Levy})}
Two non-degenerate pencils of conics are equivalent \emph{modulo}
$PGL(3,\R)$ if and only if they have the same numbers of real and
imaginary base points of each multiplicity. 
\end{theorem}

The space of couples of distinct real conics
is an algebraic fiber bundle over the variety of pencils, which is a grassmannian of
the $\RP^1$'s in a $\RP^5$. The fibers are isomorphic to the space of
couples of distinct points in $\RP^1$. The
sets of couples of distinct conics with given numbers of real and
imaginary intersections of each multiplicity are, after Theorem
\ref{pencils:base_points}, exactly the inverse images of the orbits of
the variety of pencils under $PGL(\R^3)$, and are thus also smooth
real algebraic submanifolds. 

We can now state the following definition.
\begin{definition}\label{definition:rigid_isotopy}
Two couples of distinct proper conics are \emph{rigidly isotopic} if
they are connected by a path in the space of couples of distinct
proper conics, along which the numbers of real
and imaginary intersections of each multiplicity don't change.
\end{definition}

We will now show that \emph{rigidly isotopic} implies \emph{ambient
  isotopic}. We first show it for some special rigid isotopies.
\begin{definition}
Let $f,g$ be two non-degenerate non-proportional quadratic forms.
We define a \emph{sliding} for $[f]$ and $[g]$ as a path of the
form $t \mapsto ([f+t kg],[g])$ (or $t \mapsto ([f],[g+t kf])$) for $t$ in
a closed interval containing $0$; no $t$ with $f+t k g$ (resp. $g+t
k f$) degenerate; and $k$ some real number.
\end{definition}

Let $\alpha$ be an homeomorphism of $\RP^2$. 
For a couple $([f],[g])$ of non-empty proper conics, we write
$\alpha([f],[g])$ for the couple of algebraic conics corresponding to $(\alpha([f=0]),\alpha([g=0]))$.
\begin{lemma}\label{lemma:riap}
Any sliding, for a couple of non-empty conics, lifts to an ambient isotopy.
\end{lemma}
For a sliding:
\[
t \mapsto ([f+t k g],[g]),\quad t\in I
\]
with $[f=0]$ and $[g=0]$ non-empty,
this means that 
there exists a family of homeomorphisms
$\beta_t$ of $\RP^2$, with $\beta_0=\id$ and $\beta_t([f],[g])=
([f+t k g],[g])$.

\begin{proof2}
Let $B$ be the set of the base points of the pencil of $[f]$ and $[g]$.
A stratification of $\RP^2 \times I$ is given by:
\begin{align*}
S_1&=B \times I\\
S_2&=\left([g=0]\times I\right) \setminus S_1\\
S_3&=
\left\lbrace
({\bf p};t) \; | \; (f+t k g)({\bf p})=0 \right\rbrace
\setminus 
S_1\\
S_4&= \left(\RP^2 \times I\right) \setminus
  \left(S_1 \cup S_2 \cup S_3\right)
\end{align*}
One checks this stratification is Whitney. 
The projection from $\RP^2\times I$ to $I$ is a proper stratified
submersion. The lemma now follows, by direct application of Thom's
isotopy lemma, as it is stated in \cite{Goresky:MacPherson}. 
\end{proof2}

\begin{lemma}\label{lemma:basic_ri}
Consider two couples of distinct proper non-empty conics.
If they are rigidly isotopic, then they can also be connected by a
rigid isotopy 
$\alpha_t(s_t)$ where
\begin{itemize}
\item $s_t$ is a sequence of slidings along one given pencil.
\item $\alpha_t$ is a path in $PGL(3,\R)$ with $\alpha_0=\id$
\end{itemize}
\end{lemma}

\begin{proof2}
Let $(C_0,D_0)$ and $(C_1,D_1)$ be the couples of conics, and
$t\mapsto (C_t,D_t)$, $t\in [0;1]$ be the rigid isotopy that connects
them.
It projects to a path in one $PGL(3,\R)$-orbit of the variety of pencils. 
This path lifts to a path $\alpha_t$ of $PGL(3,\R)$ with
$\alpha_0=\id$ (indeed, the group is a principal fiber bundle over
each orbit; specially, it is a locally trivial fiber bundle:
\cite{Brocker:tomDieck}, ch. I, 4).

The mapping 
$t \mapsto \alpha_t^{-1}(C_t,D_t)$ 
is a rigid isotopy drawn inside one pencil of conics. Such an isotopy is easy to
describe: the pencil is a space $\RP^1$ with a finite set $\Gamma$ of
degenerate conics. Let $E=\RP^1\setminus \Gamma$. A rigid isotopy inside
the pencil is exactly a path in $E
\times E \setminus \Diag(E \times E)$.
There exists a finite sequence $s_t$ of horizontal and vertical paths,
\emph{i.e.} of slidings, having also origin
$(C_0,D_0)$ and extremity 
$\alpha_1^{-1}(C_1,D_1)$.

Consider now $\alpha_t(s_t)$. This is a
rigid isotopy connecting $(C_0,D_0)$ to $(C_1,D_1)$.
\end{proof2}

\begin{theorem}\label{rigid:ambiant}
Two couples of distinct proper non-empty conics that are rigidly isotopic are also
ambient isotopic.
\end{theorem}
\begin{proof2}
Let $(C_0,D_0)$ and $(C_1,D_1)$ be rigidly isotopic. 
Consider a path $t \in [0,1] \mapsto \alpha_t(s_t)$ connecting them, as in 
Lemma \ref{lemma:basic_ri}.
After Lemma \ref{lemma:riap}, $s_t$ lifts to an ambient isotopy
$\beta_t$ with $\beta_0=\id$. Then $\alpha_t \circ \beta_t$ is an ambient isotopy
carrying $(C_0,D_0)$ to $(C_1,D_1)$.
\end{proof2}

\subsection{Orbits of pencils of conics}

\begin{figure}
\centering
\begin{tabular}{c@{\qquad}c}
      \includegraphics[scale=0.48]{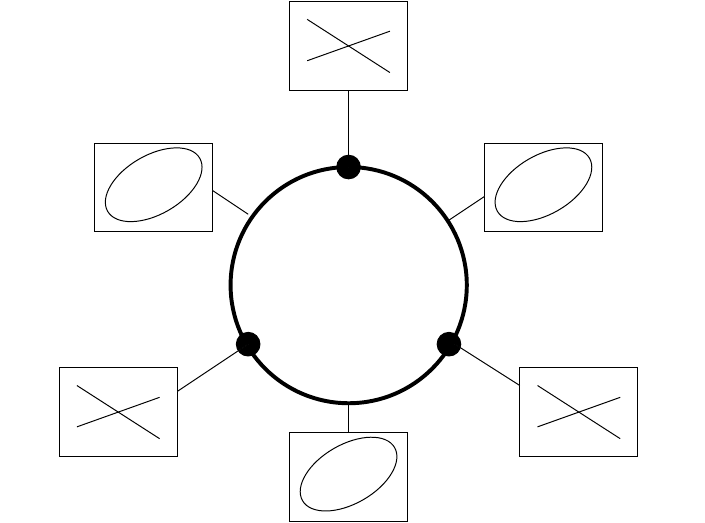} & 
      \includegraphics[scale=0.48]{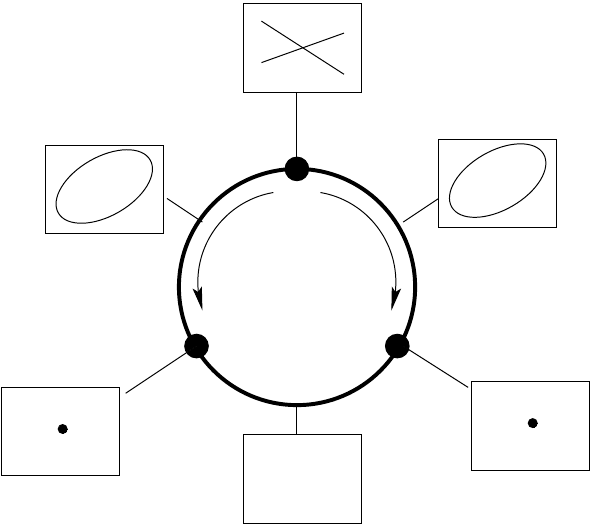} \\
{\PI}& {\PIa} \\
\multicolumn{2}{c}{
\begin{tabular}{c}
      \includegraphics[scale=0.48]{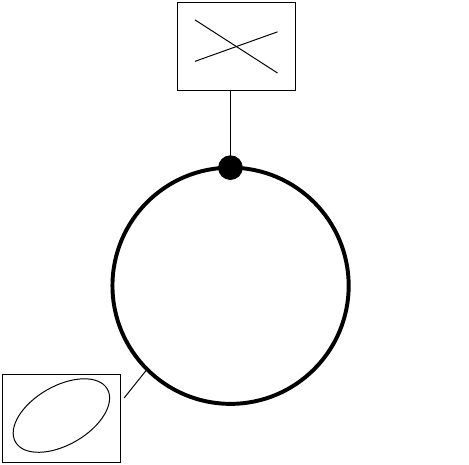} \\
      {\PIb}
\end{tabular}
}\\
      \includegraphics[scale=0.48]{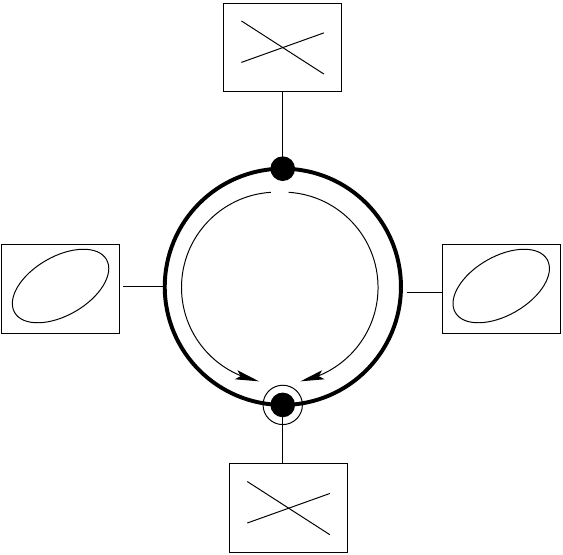} & 
      \includegraphics[scale=0.48]{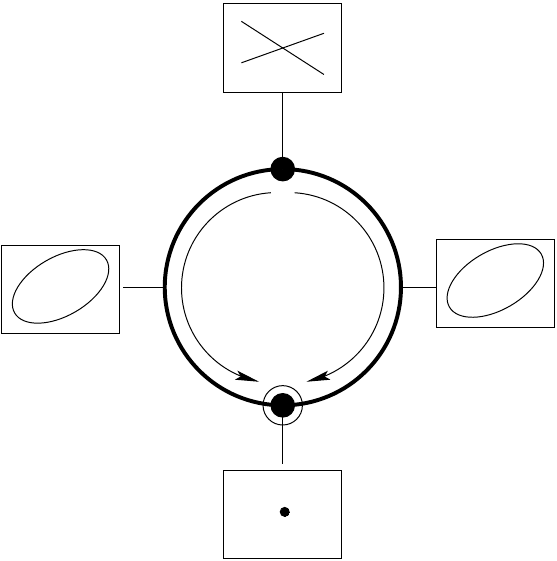} \\
{\PII}& {\PIIa} 
\end{tabular}
\caption{Pencils of conics up to projective equivalence (beginning).}\label{pencils:partI}
\end{figure}

\begin{figure}[htpb]
\begin{tabular}{c@{\qquad}c}
      \includegraphics[scale=0.48]{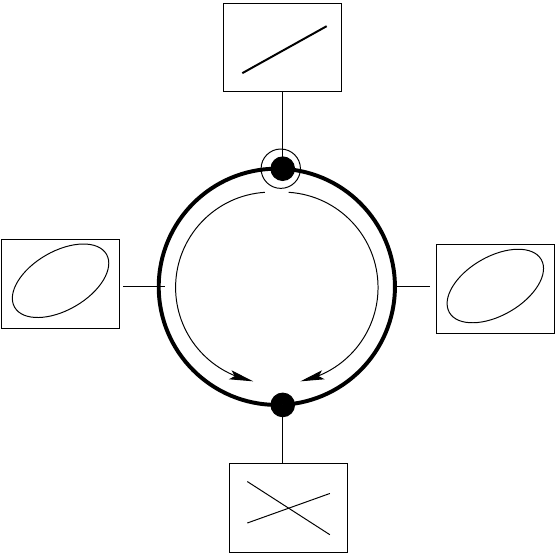} & 
      \includegraphics[scale=0.48]{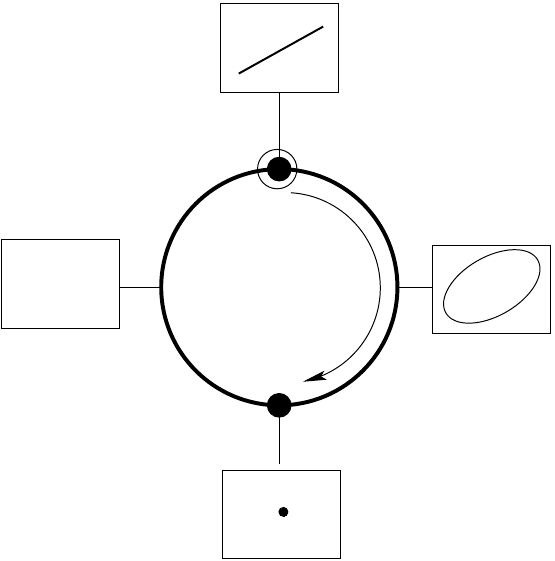} \\
      {\PIII} & {\PIIIa} \\
      &\\
      \includegraphics[scale=0.48]{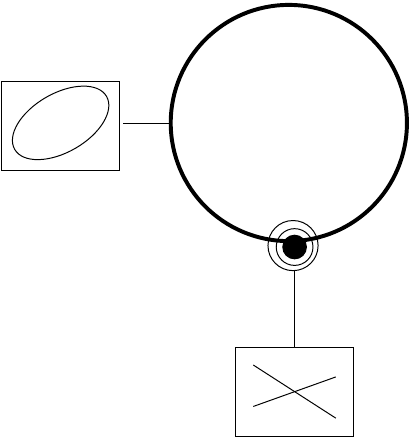} & 
      \includegraphics[scale=0.48]{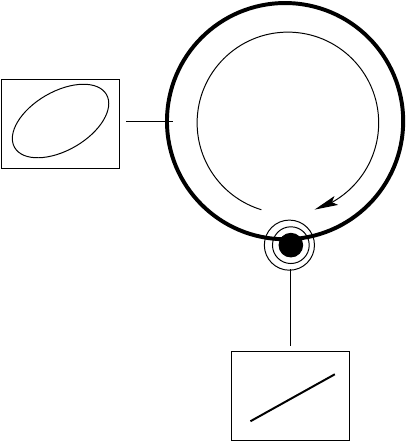}\\
      {\PIV} & {\PV} 
\end{tabular}
\caption{Pencils of conics up to projective equivalence (end).}\label{pencils:partII}
\end{figure}

After \cite{Degtyarev,Levy}, there are nine orbits of non-degenerate
pencils of conics under the action of
$PGL(3,\R)$. We follow Levy's nomenclature \cite{Levy} for them. It is
presented in the following table, where the second and third lines
display the multiplicities of the real and imaginary base
points. For instance, $211$ stands for one base point of multiplicity $2$ and
two base points of multiplicity $1$.
\[
\aaeccarray{|c|c|c|c|c|c|c|c|c|c|}
{\text{\small Orbit} &
{\PI} & {{\PIa}} & {\PIb} & {\PII}& {\PIIa} & {\PIII} & {\PIIIa} & {\PIV} & {\PV}}
{
\text{\small real\ points} &
1111 & - & 11 & 211 & 2 & 22 & - & 31 & 4\\
\hline
\text{\small imaginary\ points} & - & 1111 & 11 & - & 11 & - & 22 & -
& - \\
}
\]
We will also use the representatives of the orbits provided by Levy
\cite{Levy}. Each representative is given by a pair of generators of the
corresponding pencil of quadratic forms. They are presented in Table
\ref{table:orbit_representatives}.
\begin{table}[t]
\[
\aaeccarray{|c||c|c|c|}
{\text{\ Orbit\ } & f_0 & g_0}
{
{\PI}& x^2-y^2 & x^2-z^2 \\
{\PIa} & x^2+y^2+z^2 & xz \\
{\PIb} & x^2+y^2-z^2 & xz \\
{\PII}&  yz & x(y-z) \\
{\PIIa} & y^2+z^2 & xz \\
{\PIII} & xz & y^2\\
{\PIIIa} & x^2+y^2 & z^2 \\
{\PIV} & xz-y^2 & xy \\
{\PV} & xz-y^2 & x^2\\
}
\]
\caption{Levy's representatives for each orbit of pencils. Each representative is the
pencil generated by $[f_0]$ and $[g_0]$.}\label{table:orbit_representatives}
\end{table}
We
provide, in figures \ref{pencils:partI} and \ref{pencils:partII}, graphical representations of characteristic features of
the pencils in each orbit. This has two goals: finding how to
discriminate between the different orbits of pencils, and determining the possible
rigid isotopy classes corresponding to each orbit of pencils. 
Each pencil is displayed as a circle, as it
is topologically. In addition, the following information is represented:
\begin{itemize}
\item the degenerate conics of the
  pencil. They are given by roots of the characteristic form, so the
  multiplicity of this root is also indicated, following the encoding
  shown in Figure \ref{multiplicities}.
\begin{figure}
\centering
\includegraphics[scale=0.75]{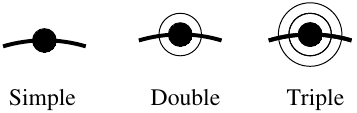}
\caption{Degenerate conics corresponding to multiple roots of the
  discriminant, in the representations of the pencils.}\label{multiplicities}
\end{figure} 
\item the nature of the proper conics (empty or non-empty) and of the degenerate conics (pair of lines, line
  or isolated point). The nature of the proper conics is constant on
  each arc between two degenerate conics.
\begin{figure}
\centering
\includegraphics[scale=0.75]{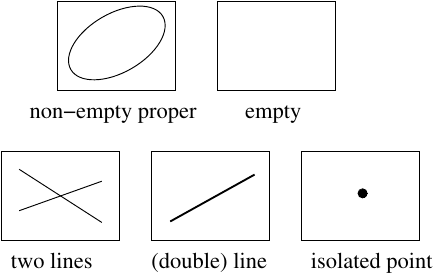}
\caption{Nature of the conics, 
in the representations of the pencils.}\label{nature}
\end{figure} 
\item In the case where the conics of one arc are nested, we indicate,
  by means of an arrow, which are
  the inner ones.
\end{itemize}
These features are conserved under projective equivalence. Thus the
representations are established by considering Levy's representative.

\subsection{Rigid isotopy classification for pairs}

We 
first classify \emph{pairs} of proper conics, that is: couples of
distinct proper conics, under rigid isotopy \emph{and permutation of
  the two conics}. Later, for each pair class, we will check whether it is also a
couple class (that is: the exchange of the two conics corresponds to
a rigid isotopy) or it splits into two couple classes.

Any pencil of conics is cut into arcs by its degenerate conics. Two proper
conics are either on a same arc, or on distinct arcs on the pencil
they generate.

\begin{lemma}\label{lemma:pairs}
If a pencil of conics has (at least) two arcs of non-empty conics, then there are two
equivalence classes for pairs of conics generating it.
They correspond to the following situations:
\begin{itemize}
\item the conics are on a same arc.
\item the conics are on distinct arcs.
\end{itemize}
If the pencil has only one arc with proper non-empty conics, there is only one class.
\end{lemma}

\begin{proof2}
The orbit of pencils is assumed to be fixed.
Because of Lemma \ref{lemma:riap}, to get (at least) one representative
for each class, it is enough:
\begin{itemize}
\item to choose arbitrarily one conic on each arc and consider all the possible
  pairs of these conics.
\item to choose arbitrarily two conics on each arc and consider these pairs for
  each arc.
\end{itemize}
But one observes that for a pencil in one of the orbits {\PIa},
{\PII}, {\PIIa}, {\PIII}, there is a projective automorphism that leaves it
globally invariant and 
exchanges its
two arcs (the two arcs bearing non-empty conics for orbit {\PIa}). 
Again, this is proved by considering only Levy's
representatives: the reflection
$x \leftrightarrow -x$ is suitable. 
Similarly, a pencil in orbit {\PI} is left globally
invariant by some projective automorphism that permutes cyclically
the three arcs. For Levy's representative, one can take the cyclic
permutation of coordinates: $x \mapsto y \mapsto z \mapsto x$.

Pencils in the four other orbits have only one arc with non-empty
proper conics. 

We have shown it is enough:
\begin{itemize}
\item to choose arbitrarily one arc with non-empty conics and two
  conics on this arc.
\item to choose arbitrarily two arcs and one conic on each arc.
\end{itemize}
This gives nine representatives for the pairs of conics on a same arc,
denoted with {{\CIN}}, \ldots, {\CVN}  ({\CN}  like \emph{neighbors}) and five 
representatives for pairs of conics on distinct arcs, denoted with {\CIS},
{\CIaS}, {\CIIS}, {\CIIaS}, {\CIIIS}  ({\CS}  like \emph{separated}).  

Now it remains to check that for orbits {\PI}, {\PIa}, {\PII} {\PIIa} and {\PIII}, the {\CS}--representative and the
{\CN}--representative are not equivalent. We use that a rigid isotopy
conserves the topological type of $(\RP^2,[f=0],[g=0])$, after Theorem
\ref{rigid:ambiant}. To distinguish between {{\CIN}}    and {\CIS}, one can count the number of
connected components of the complement of $[f=0] \cup [g=0]$: there
are $6$ in the first case and $5$ in the second, the topological types
are different, so are the rigid isotopy classes. For the other four
orbits of pencils, one conic lies in the inside of the other (at least
at the neighborhood of the double point for {\PII}) for the
{\CN}-representative, while there is no such inclusion for the {\CS}-representative.
\end{proof2}

\begin{corollary}
There are $14$ equivalence classes for pairs of proper
non-empty conics under rigid isotopy and exchange. Representatives for
them are given in Table \ref{table:r_isotopy}. They correspond to the
graphical representations displayed in Figures \ref{fig:generic} and \ref{fig:singular}.
\end{corollary}

\begin{table}
\[
\aaeccarray{|c|c|c|}
{
\text{\ class\ } & f & g
}{
{{\CIN}}               & 3x^2-2y^2-z^2     & 3x^2-y^2-2z^2          \\
{\CIS}               & 3x^2-2y^2-z^2     &  x^2-2 y^2+z^2         \\
{\CIaN} (*)          & x^2+y^2+z^2+3xz   & x^2+y^2+z^2+4xz        \\
{\CIaS}              & x^2+y^2+z^2+3xz   &  x^2+y^2+z^2-3xz       \\
{\CIbN}              & x^2+y^2-z^2+xz    & x^2+y^2-z^2-xz         \\
{\CIIN} (*)          & yz+xy-xz          &  yz+2xy-2xz            \\
{\CIIS}              & yz+xy-xz          &  yz-xy+xz              \\
{\CIIaN} (*)         & y^2+z^2+xz      &  y^2+z^2+2 xz            \\
{\CIIaS}             & y^2+z^2+xz      &  y^2+z^2-xz           \\
{\CIIIN} (*)         & xz+^2          &  xz+2 y^2               \\
{\CIIIS}             & xz+y^2          &  xz-y^2              \\
{\CIIIaN} (*)        & x^2+y^2-z^2       &  x^2+y^2-2 z^2         \\
{\CIVN}              & xz-y^2+xy         &  xz-y^2-2 xy           \\
{\CVN} (*)           & xz-y^2-x^2        &  xz-y^2+x^2            \\
}
\]
\caption{The rigid isotopy classes.}\label{table:r_isotopy}
\end{table}

\subsection{Rigid isotopy classification for couples}

We now derive from our classification for pairs of conics the
classification for couples of conics.

\begin{lemma}
For each of the following representatives: {{\CIN}}, {\CIS}, {\CIaS},
{\CIbN}, {\CIIS}, {\CIIaS}, {\CIIIS}, {\CIVN}, there is a rigid isotopy that swaps the two conics. As a consequence, each of
these classes for pairs is also a class for couples.
\end{lemma}

\begin{proof2}
For {{\CIN}}, {\CIS}, {\CIaS}, {\CIIS}, {\CIIaS}, {\CIIIS}, it is enough to exhibit
projective automorphisms that stabilize the corresponding Levy's
representative and swap two arcs of non-empty proper conics. It was
already done in the proof of Lemma \ref{lemma:pairs}, except for {{\CIN}}
and {\CIS}. For them, the reflection $y \leftrightarrow z$ is
convenient. 

For {\CIbN} and {\CIVN}, it is enough to exhibit projective automorphisms
that stabilize the corresponding Levy's representative and reverse the
pencil's orientation. For {\CIbN}, the reflection $x \mapsto -x$ is
convenient; for {\CIVN}, one may use the transformation $x\mapsto -x, z \mapsto -z$.
\end{proof2}

\begin{lemma}\ 
\begin{itemize}
\item Each of the classes of pairs {\CIaN}, {\CIIaN}, {\CIIIN}, {\CIIIaN}, {\CVN}
splits into two classes for couples, corresponding to one conic lying
inside the other (except for the base points).
\item The class of pairs {\CIIN} also splits into two classes for couples,
  corresponding to one conic lying inside the other in a neighborhood
  of the double point (except the double point itself).
\end{itemize} 
\end{lemma}

\begin{proof2}
The property that one conic lies inside the second is conserved under ambient homeomorphism, and thus under rigid
isotopy. The same holds for inclusion at the neighborhood of a double
intersection point. 

Thus it is enough to consider the representatives of the given pair
classes and check the inclusion to show the theorem. The computations
are trivial, hence we omit them.
\end{proof2}

\begin{theorem}
There are $20$ classes of couples under rigid isotopy. A set of
representatives is given by Table \ref{table:r_isotopy}, where the
reader should add the couple obtained by swapping $f$ and $g$ for each of the
lines marked with $(*)$.  
\end{theorem}

\begin{corollary}
The ambient isotopy classes for couples of conics are the following
unions of rigid isotopy classes: 
\begin{itemize}
\item classes where the two conics can be swapped: {{\CIN}}, {\CIS}, {\CIaS},
${\CIbN} \cup {\CIVN}$, {\CIIS}, {\CIIaS} and {\CIIIS}.
\item pair classes splitting into two classes for couples, one with
  $[f=0]$ inside $[g=0]$, one with $[g=0]$ inside $[f=0]$: 
${\CIaN} \cup {\CIIIaN}$, {\CIIN}, ${\CIIaN} \cup {\CVN}$ and  {\CIIIN}.  
\end{itemize}
\end{corollary}

\begin{proof2} \emph{(sketch)} 
One shows that {\CIbN} and {\CIVN} are ambient isotopic by building
explicitly a
homeomorphism\footnote{This is enough. Indeed, the set of the
  homeomorphisms of $\RP^2$ is connected, so any homeomorphism is the
  extremity of some ambient isotopy.}
 of $\RP^2$ sending a representative of the first to a
representative of the second.
Details on how to do it are tedious, we skip them.
\emph{Idem} (with $[f=0]$
inside $[g=0]$) for {\CIaN}
and {\CIIIaN}, and for {\CIIaN} and {\CVN}. 

Next, one shows the displayed rigid isotopy classes or unions of rigid
isotopy classes are not equivalent \emph{modulo}
ambient isotopy. This is done by considering topological invariants of the triples
$(\RP^2, [f=0],[g=0])$ which take different values on the  $15$
representatives. 

Let $C$ be the conic $[f=0]$
(resp. $D$ the conic $[g=0]$), $I$ (resp. $J$) its inside and
$\bar{I}$ (resp. $\bar{J}$) the topological closure of this inside.

Then one checks that the numbers of connected components of the four 
following sets are suitable for separating the ambient isotopy
classes: 
\[
C \cap D,\quad
 \RP^2 \setminus (C \cup D), 
\quad
I \setminus \bar{J}, \quad
J \setminus \bar{I} 
\]
\end{proof2}

\begin{remark}
It is legitimate curiosity to compare these isotopy
  classes of couples of conics with the isotopy classes of projective
  quartic curves 
  presented in \cite{Korchagin:Weinberg},
  the union of two conics being a quartic. One then observes that 
{{\CIN}} corresponds to $17p$, 
{\CIS} to $16p$, 
{\CIaS} to $22p$, 
{\CIIS} to $34p$, 
{\CIIaS} to $44p$, 
{\CIIIS} to $38p$, 
${\CIaN} \cup {\CIIIaN}$ to $21p$, 
{\CIIN} to $36p$ 
and 
${\CIIaN} \cup {\CVN}$ to $43p$. 
Finally both ${\CIbN} \cup \CIVN$ and {\CIIIN} correspond to $18p$.
\end{remark}


\section{Characterizing the isotopy classes by equations, inequations
  and inequalities}

\subsection{Preliminaries}

\subsubsection{The invariants and covariants of two ternary quadratic forms}\label{subsection:invariants}

Invariants and covariants (see \cite{Olver} for a modern reference
about classical invariant theory)
are the convenient objects to discriminate, by means of equations and
inequalities, between the different orbits for couples of complex
conics under the group $PGL(3,\C)$. 

Invariants and covariants of a couple of quadratic ternary forms have been calculated by the classics, and can be found in Glenn's
book \cite{Glenn} or Casey's treatise \cite{Casey}. 

\begin{proposition}
The algebra of invariants of a couple of ternary quadratic forms is
freely generated by the coefficients of the characteristic
form\footnote{The analogue assertion is still true for a couple of
  quadratic forms in $n$ variables, for any $n$.}.
\end{proposition}

The invariants alone are not sufficient to discriminate between the
complex orbits. One has to consider the covariants. Some
remarkable
covariants of a couple of ternary quadratic forms are:
\begin{itemize}
\item The \emph{apolar covariant of the tangential quadratic forms}
  $\tilde{f}$ and
  $\tilde{g}$. We will denote it with $F$. This is a quadratic form that
  depends quadratically on $f$, as on $g$.
\item The \emph{autopolar triangle covariant} $G$, a cubic form that
  is also cubic in $f$, as in $g$, and that always factorizes as a
  product of three linear forms. When $[f=0]$ and $[g=0]$ have four
  distinct intersections, they are equations of the sides of the
  unique autopolar triangle associated to them (see \cite{Berger},
  14.5.4 and 16.4.10). 
\end{itemize}

\begin{proposition}
The algebra of covariants of a couple of ternary quadratic forms $(f,g)$ is
generated by the invariants, the ground forms $f$ and $g$, the apolar
covariant $F$ and the autopolar triangle covariant $G$.
\end{proposition}

The covariant $F$ will not be needed in this paper, but $G$ will be
used. We now 
explain how to derive a formula for it.
Consider a generic
couple of forms $f,g$. 
Let $t_1,t_2,t_3$ be the three roots of $\Disc(t f+g)$.
Each of the $t_i f+g$ has rank two. Their respective associated tangential
quadratic forms have all 
rank one: they are the squares of three linear forms $p_1,p_2,p_3$ of the dual
space, and the associated points $[p_1],[p_2],[p_3]$ of $\P(\R^3)$ are
exactly the
vertices of the autopolar triangle.
The sides of the triangle are obtained as the zero loci of the
product of determinants:
\[
\det(p_1,p_2,p)\det(p_1,p_3,p)\det(p_2,p_3,p).
\]
Working in coordinates, with 
$p_i=p_{i1}X+p_{i2}Y+p_{i3}Z$, where $X,Y,Z$ are coordinates on
${\R^3}^*$ dual to $x,y,z$, 
one expands this and replace the products $p_{ij} p_{ik}$
by the corresponding term given by the equality
\[
p_i^2=\widetilde{f t_i+g}.
\]
This product is antisymmetric in $t_1,t_2,t_3$, and thus can be
divided by the Vandermonde determinant
$(t_1-t_2)(t_1-t_3)(t_2-t_3)$. The quotient happens to be free of
$t_i$'s: it is, up to a rational number in factor, the covariant $G$.

One finds that the formula for this covariant can be displayed shortly.
Denote:
\[
\widetilde{t f+g}
=
\tilde{f} t^2
+
\Omega(f,g) t
+
\tilde{g}
\]
and
\[
\tilde{f}=
\sum \tilde{a}_{ijk} X^i Y^j Z^k,\quad
\Omega(f,g)=
\sum \omega_{ijk} X^i Y^j Z^k,\quad
\tilde{g}=
\sum \tilde{b}_{ijk} X^i Y^j Z^k.
\]
Consider the matrix of their coefficients:
\[
M=
\left[
\begin{matrix}
\tilde{a}_{200}  & \tilde{a}_{020} & \tilde{a}_{002} & \tilde{a}_{011}
& \tilde{a}_{101} & \tilde{a}_{110} \\
\omega_{200}  & \omega_{020} & \omega_{002} & \omega_{011}
& \omega_{101} & \omega_{110} \\
\tilde{b}_{200}  & \tilde{b}_{020} & \tilde{b}_{002} & \tilde{b}_{011}
& \tilde{b}_{101} & \tilde{b}_{110} \\
\end{matrix}
\right]
\]
and label its columns with $1,2,3,\bar{1},\bar{2},\bar{3}$. Label the maximal minor 
corresponding to columns $i,j,k$ with $[ijk]$.
Then the autopolar triangle covariant is, up to a rational number in factor, 
\begin{multline*}
- [\bar{1}\,2\,3] x^3
- [1\,\bar{2}\,3] y^3
- [1\,2\,\bar{3}] z^3+ ([1\,\bar{1}\,3]+2[\bar{3}\,\bar{2}\,3])x y^2 \\
+ ([1\,2\,\bar{1}]+2[\bar{2}\,2\,\bar{3}])x z^2 
+ ([\bar{2}\,2\,3]+2[\bar{1}\,\bar{3}\,3])y x^2 
+ ([1\,2\,\bar{2}]+2[1\,\bar{1}\,\bar{3}])y z^2\\ 
+ ([\bar{3}\,2\,3]+2[\bar{1}\,2\,\bar{2}])z x^2 
+ ([1\,\bar{3}\,3]+2[1\,\bar{2}\,\bar{1}])z y^2 
+ ([1\,2\,3]+4[\bar{1}\,\bar{2}\,\bar{3}]) xyz.
\end{multline*}

\subsubsection{Resultants}

Let $U(t)$ and $V(t)$ be two univariate polynomials. 
Remember that their resultant $\Res(U,V)$ is the determinant of
their Sylvester matrix: the matrix of the coefficients of degree
$\deg(U)+\deg(V)-1$ down to $0$ of 
\[t^{\deg(V)-1}U,t^{\deg(V)-2}U,
\ldots, U, t^{\deg(U)-1}V,t^{\deg(U)-2}V,
\ldots, V.
\]

A few classical formulas about resultants will be needed.
\begin{lemma}\label{res:exchange}
One has
\[
\Res(U,V)= (-1)^{\deg(U)\deg(V)} \Res(V,U).
\]
\end{lemma}

\begin{lemma}\label{res:roots}
Let $c$ be the leading coefficient of $U$. Then
\[
\Res(U,V)=c^{\deg(V)} \prod V(\rho)
\]
where the product is carried over the complex roots $\rho$ of $U$, counted
with multiplicities\footnote{\emph{e.-g.}
  a double
  real root should be here counted as two roots.}.
\end{lemma}

\begin{lemma}\label{res:product}
Let $U(t),V(t),W(t)$ be three univariate polynomials. Then
\[
\Res(U,VW)=\Res(U,V)\Res(U,W).
\]
\end{lemma}

And last:
\begin{lemma}\label{res:remainder}
Let $U(t),V(t)$ be two univariate polynomials, and $W$ the remainder in the euclidean
division of $U$ by $V$. Let $c$ be the leading coefficient of $V$. Then
\[
\Res(U,V)=(-1)^{\deg(U)\deg(V)} c^{\deg(U)-\deg(W)} \Res(V,W).
\]
\end{lemma}
This is Lemma 4.27 in \cite{Basu:Pollack:Roy}, where a proof is provided.

\subsubsection{Descartes' law of signs}

Let $U(t)$ be an univariate polynomial. 
Then Descartes' law of signs give some insight about its number $\mathcal{N}(U)$ of
positive real roots, counted with multiplicities.

Consider the sequence of the signs ($+$'s and $-$'s) of the (non-zero)
coefficients of $U$ and denote with $\mathcal{V}(U)$ the number of changes in
consecutive terms. 
The following lemma is Descartes' Law of signs. It can be found as Theorem
2.34 in \cite{Basu:Pollack:Roy}.
\begin{lemma}
One has
$\mathcal{V}(U) \geq \mathcal{N}(U)$, and $\mathcal{V}(U)-\mathcal{N}(U)$ is even.
\end{lemma}

Only the following particular consequence will be needed in the sequel:
\begin{lemma}\label{Descartes:hyperbolic}
Let $U(t)=u_3 t^3 + u_2 t^2 + u_1 t + u_0$ of degree $3$. Suppose $U$ has all its
roots real and non-zero. Then they have all the same sign if and only if: $u_3 u_1 >0$ and $u_2
u_0 >0$.
\end{lemma}
It is obtained by applying Descartes' law of signs to $U(t)$ and $U(-t)$.

\subsubsection{Subresultant sequences}\label{subsection:subresultants}

Here we briefly introduce another tool: subresultant sequences. More
details about them can be found
in the book \cite{Basu:Pollack:Roy}. 

Let $U(t),V(t)$ be two univariate polynomials. 
One wants to know on how many\footnote{When dealing with subresultant
  sequences, the multiplicities of the roots are not taken into
  account, \emph{e.-g.} a double root will be counted as one root.} of the (real) roots of $V$ the
polynomial $U$ is positive, negative, and zero. 
The \emph{Sturm query of $U$ for $V$} is defined as the number of
roots of $V$ making $U>0$, \emph{minus} the number of roots of $V$ making
$U<0$. This information is easily accessible once one knows the signs
of the $\deg(V)+1$ \emph{signed subresultants principal coefficients} of $V$ and
$W$, where $W$ is the remainder in the euclidean division of $U \cdot V'$ by $V$. 

We give the formulas for these signed subresultant principal
coefficients, and the procedure for getting the Sturm query from their
signs, only for the particular case needed: when $V$ has degree $3$. 
Write
\[
\begin{array}{cll}
V &=v_3 t^3 + & v_2 t^2 + v_1 t + v_0\\
W &=          & w_2 t^2 + w_1 t + w_0.
\end{array}
\]
Then
\[
\begin{array}{c@{\qquad}c}
\sr_3(V,W)=v_3,&
\sr_2(V,W)=w_2,\vspace{5mm}\\
\sr_1(V,W)=
\left\vert
\begin{matrix}
v_3 & v_2 & v_1 \\
0   & w_2 & w_1 \\
w_2 & w_1 & w_0
\end{matrix}
\right\vert,&
\sr_0(V,W)=
\left\vert
\begin{matrix}
v_3 & v_2 & v_1 & v_0 & 0 \\
0   & v_3 & v_2 & v_1 & v_0  \\
0   & 0   & w_2 & w_1 & w_0 \\
0   & w_2 & w_1 & w_0 & 0 \\
w_2 & w_1 & w_0 & 0   & 0
\end{matrix}
\right\vert.
\end{array}
\]
Note that $\sr_0(V,W)=-\Res(V,W)$, the opposite of the resultant of $V$
and $W$.
The Sturm query is obtained from the sequence of the signs of
$\sr_3,\sr_2,\sr_1,\sr_0$ the following way\footnote{for this specific
case with $4$ terms in the sign sequence.}:
\begin{enumerate}
\item If there is a pair of consecutive zeros, remove it and change
  the signs that were following to their opposites.
\item From the resulting sequences of consecutive non-zero terms, compute the difference:
  number of sign permanences (identical consecutive terms,
  $++$ or $--$) \emph{minus} number of sign exchanges (opposite consecutive terms, $+-$ or $-+$). 
This gives the Sturm query\footnote{So for instance, the sign sequence $+0-0$ has no sign permanence, nor
sign change (because there are no consecutive non-zero terms). For the
sign sequence $+00-$, the Sturm query is computed as for $++$: one
permanence, no change, this gives $1$.}.
\end{enumerate}

\subsection{Discriminating between the orbits of pencils}

In this section, we give the equations and inequations characterizing
the couples $(f,g)$ of non-degenerate quadratic
forms generating a pencil of each of the orbits.

We first use invariants and covariants whose vanishing depends only of
the generated pencil, that is those $C$ that, besides the good behavior
with respect to the action of $SL(3,\C)$:
\[
C(f \circ \theta, g \circ \theta; x,y,z;t,u)=C(f,g;\theta(x,y,z);t,u) \qquad
\forall \theta \in SL(3,\C)
\] 
are covariant with respect to combinations of $f$ and $g$:
\[
C(\theta(f,g); x,y,z;\theta(t,u))=C(f,g;x,y,z;t,u) \qquad \forall
\theta \in SL(2,\C).
\]
Such objects are called \emph{combinants}.

Obviously, the characteristic form $\Phi$ and its covariants are combinants.
Remember that the algebra of the covariants of a binary cubic form
$\Phi(t,u)$ is generated by the ground form $\Phi$, its
discriminant\footnote{The discriminant of the characteristic form,
  $\Disc(\Phi)$, is called the \emph{Tact invariant} by the classics, because it
  vanishes exactly when the two conics are tangent \cite{Casey}.},
and its Hessian determinant, which are 
\[
\Disc(\Phi)=
\frac{\operatorname{Res}(\phi,\phi')}{27 \Phi_{30}},
\quad 
H(t,u)
=
\left\vert
\begin{matrix}
\frac{\partial^2 \Phi}{dt^2}  & \frac{\partial^2 \Phi}{dt\, du}   \\
\frac{\partial^2 \Phi}{dt \, du}  & \frac{\partial^2 \Phi}{d u^2}   
\end{matrix}
\right\vert
\]
(the division is a simplification in the definition of the
discriminant, that is: there remains no $\Phi_{30}$ at the denominator).
The covariant $G$ is also a combinant.

The vanishing or non-vanishing of each of the combinants are
properties of the orbits of pencils of conics. 
The sign of $\Disc(\Phi)$ is also invariant on each orbit of pencils
of conics (because $\Disc(\Phi)$ has even degree in $f$ as well
as in $g$).
Thus we just evaluate the combinants on Levy's representatives, and we get the
following result:
\begin{proposition}
Let $f,g$ be two non-proportional non-degenerate ternary quadratic forms. 
\begin{itemize}
\item If $\Disc(\Phi)<0$ then $f,g$ generate a pencil in orbit {\PI} or
  {\PIa}.
\item If $\Disc(\Phi)>0$ then $f,g$ generate a pencil in orbit {\PIb}.
\item If $\Disc(\Phi)=0$ then $f,g$ generate a pencil in one of the
  six other orbits. The following table indicates how the vanishings of
  $H$ and $G$ discriminate further between the orbits of pencils
  (under the hypothesis that the discriminant vanishes):  
\[
\begin{array}{|c|c|c|}
\cline{2-3}
\multicolumn{1}{c|}{} & H \neq 0 & H=0 \\
\hline
G\neq 0             & {\PII},{\PIIa}          & {\PIV}         \\
\hline
G=0                 & {\PIII},{\PIIIa}        & {\PV}          \\
\hline
\end{array}.
\]
\end{itemize}
\end{proposition}

\begin{remark}
The fact that the coefficients of $G$ are linear combinations of
maximal minors of the matrix $M$ defined in \ref{subsection:invariants} suggests
that the vanishing of $G$ is equivalent to: \emph{$M$ takes rank two}. This is true. To see this, consider the image of the (complex) pencil
generated by $[f]$ and $[g]$ by the quadratic mapping ``tangential
quadratic form'' from $\P({S^2 \C^3}^*)$ to $\P(S^2 \C^3)$. It is an
irreducible conic, thus either a proper conic or a line. One checks
on Levy's representative that it is a line exactly when $G=0$ (see also
\cite{Berger}, 16.5.6.2). Finally, remark that the rows of $M$
are the coordinates of generators
of the linear span of this conic.
\end{remark}

It remains now to discriminate between {\PI} and {\PIa}, between
{\PII} and
{\PIIa} and between {\PIII} and {\PIIIa}.
For this we use that the numbers of degenerate conics of each
type (pair of lines, isolated point or double line) in a pencil
characterize its orbit, as shown in the table below, established by
considering figures \ref{pencils:partI} and \ref{pencils:partII}.
\[
\begin{tabular}{|c||c|c||c|c||c|c|}
\hline
orbit of pencils         &{\PI}&{\PIa}&{\PII}&{\PIIa}&{\PIII}&{\PIIIa}\\[3pt]
\hline
num. pairs of line   & 3 & 1  & 2  & 1   &  1  & 0 \\   
\hline
num. isolated points & 0 & 2  & 0  & 1   &  0  & 1 \\
\hline
num. (double) lines  & 0 & 0  & 0  & 0   &  1  & 1 \\
\hline
\end{tabular}
\]

This has an algebraic translation.
Consider
\[
\det(v \cdot I - \operatorname{Matrix}(t f + g))
\]
that expands into
\[
v^3 - \mu(t) v^2 + \psi(t) v - \phi(t). 
\]
A degenerate conic of the pencil corresponds to a parameter $t$
that annihilates $\phi$, and is
\begin{itemize}
\item a pair of lines when
  $\operatorname{Matrix}(t f + g)$ has one eigenvalue positive, one
  negative, and one zero. Then $\psi(t)<0$.
\item an isolated point when the matrix has an eigenvalue zero and the two
  other both positive or both negative. Then $\psi(t)>0$.
\item a single line when the matrix has two
  eigenvalues zero, and one non-zero. Then $\psi(t)=0$.
\end{itemize}
Thus the discriminations can be performed by a Sturm query of $\psi$
for $\phi$.
In order not to introduce denominators, we consider the euclidean
division of $\Phi_{30} \psi \phi'$ by $\phi$, instead of the
division of $\psi \phi'$ by $\phi$ suggested by \ref{subsection:subresultants}.
Set
\begin{align*}
P&=\operatorname{Remainder}(\Phi_{30} \cdot \Psi \cdot \phi',\phi)\\
 &=p_2 t^2 + p_1 t + p_0
\end{align*}
and 
\[
A_i=\sr_i(\phi,P)
\]
for $i$ between $0$ and $3$.
The consideration of the sign permanences and sign exchanges in $\Phi_{30}=A_3,A_2,A_1,A_0$ gives the
Sturm query of $\Phi_{30} \psi$ for $\phi$. 
The Sturm query of $\psi$ for $\phi$ is the same as the Sturm query of
$\Phi_{30}^2 \psi$ for $\phi$. Using that $\sr_i(\phi,\Phi_{30}
P)=\Phi_{30}^{4-i}\sr_i(\phi,P)$, we get that this Sturm query is 
obtained by considering
the sign permanences and sign exchanges in 
$\Phi_{30},\Phi_{30}A_2,A_1,\Phi_{30}A_0$; or, simpler, those in $1,A_2,\Phi_{30}A_1,A_0$.

The polynomial $A_1$ is
\[
\left\vert
\begin{matrix}
\Phi_{30} &  \Phi_{21} &  \Phi_{12} \\
0           &    p_2        &    p_1        \\
p_2         &    p_1        &    p_0   
\end{matrix}
\right\vert.
\]
And 
\[
A_0=-\operatorname{Res}(\phi,P).
\]
This simplifies. Applying Lemma \ref{res:remainder}, one gets
\[
\operatorname{Res}(\Phi_{30} \psi \phi', \phi)
=
\Phi_{30}^2 \operatorname{Res}(\phi,P).
\]
And, on the other hand, from Lemma \ref{res:product} and the definition of $\Disc(\Phi)$:
\[
\operatorname{Res}(\Phi_{30} \psi \phi', \phi)
=
 \Phi_{30}^4 
\cdot
\operatorname{Res}(\psi, \phi)
\cdot
\Disc(\Phi).
\]
Thus 
\[
A_0=-  \Phi_{30}^2 \Res(\psi,\phi) \Disc(\Phi).
\]

\subsubsection{Discriminating between {\PI} and {\PIa}}

Suppose $f$ and $g$ generate a pencil in orbit {\PI} or {\PIa}.
The Sturm query of $\psi$ for $\phi$ is 
$-3$ for orbit {\PI}
and $1$ for orbit {\PIa}. 

The assumption that $f,g$ generate a pencil in orbit {\PI} or {\PIa}
gives more information:
\begin{lemma}
If $f,g$ generate a pencil in orbit {\PI} or {\PIa}, then $A_0<0$.
\end{lemma}

\begin{proof2}
We had established that
$A_0=- \Phi_{30}^2 \Res(\psi,\phi) \Disc(\Phi)$.
For orbit {\PI} or {\PIa}, one has $\Disc(\Phi)<0$.
Moreover,
$\operatorname{Res}(\psi,\phi)<0$ because, from lemmas
\ref{res:exchange} and \ref{res:roots},  
\[
\operatorname{Res}(\psi,\phi)=
\operatorname{Res}(\phi,\psi)=
\Phi_{30}^2 
\prod \psi(\rho),
\]
where the product is carried over the three roots $\rho$ of
$\phi$. They make either $\psi$ three times negative (orbit {\PI}), either
one time negative and two times positive (orbit {\PIa}). In both cases,
the product is negative.
\end{proof2}

There is only one sign sequence giving Sturm query $-3$ and beginning
with $+$ and finishing with $-$, that is $+-+-$.
There are several sign sequences giving Sturm query $1$, beginning
with $+$, finishing with $-$:
\[
+++-\quad ++-- \quad ++0- \quad +--- \quad +0-- \quad +00-.
\]
We deduce from this the criterion stated in the following proposition.

\begin{proposition}
Let $f,g$ be non-degenerate quadratic forms generating a pencil
in orbit {\PI} or {\PIa}. 
\begin{itemize}
\item if it is orbit {\PI} then $p_2 < 0$ and $\Phi_{30}A_1 > 0$.
\item if it is orbit {\PIa}, then $p_2 > 0$, or $\Phi_{30} A_1 < 0$, or $p_2=A_1=0$.
\end{itemize}
\end{proposition}

\subsubsection{Discriminating between {\PII} and {\PIIa}}

Suppose $f$ and $g$ generate a pencil in orbit {\PII} or {\PIIa}. Note
first that $A_0=0$.
The Sturm query of $\psi$ for $\phi$ is $-2$ for orbit {\PII} and $0$ for orbit {\PIIa}. 

There is only one sign sequence with beginning with $+$, finishing
with $0$ that gives Sturm query $-2$, this is $+-+0$. 
Those giving Sturm query $0$ are 
\[
++-0 \quad 
+--0 \quad
+0+0 \quad
+0-0 \quad
+000.
\]

\begin{proposition}
Let $f,g$ be non-degenerate quadratic forms generating a pencil
in orbit {\PII} or {\PIIa}.
\begin{itemize}
\item if it is in orbit {\PII} then $p_2<0$ and $\Phi_{30} A_1 >0$.
\item if it is in orbit {\PIIa}, then $p_2=0$ or $\Phi_{30} A_1 <0$.
\end{itemize}
\end{proposition}

\subsubsection{Discriminating between {\PIII} and {\PIIIa}}

Suppose $f$ and $g$ generate a pencil in orbit {\PIII} or
{\PIIIa}. Once again, $A_0=0$.
The Sturm query of $\psi$ for $\phi$ is 
$-1$ for orbit {\PIII},
$1$ for orbit {\PIIIa}. 

The only sign sequence (beginning with $+$, terminating with $0$) giving Sturm query $-1$ is $+-00$. 
There is also only one giving Sturm query $1$, that is $++00$.

\begin{proposition}
Let $f,g$ be non-degenerate quadratic forms generating a pencil in orbit {\PIII} or {\PIIIa}.
\begin{itemize}
\item if it is in orbit {\PIII}, then $p_2 < 0$.
\item if it is in orbit {\PIIIa}, then $p_2 > 0$.
\end{itemize}
\end{proposition}

\subsection{Characterizing the rigid isotopy classes for pairs inside each pencil}

Given $f,g$ two non-proportional
non-degenerate quadratic forms, we suppose we know the orbit of the
pencil they generate.
After Lemma \ref{lemma:pairs}, one decides to which class belongs
$\{[f],[g]\}$ by looking whether or not $[f]$ and $[g]$ are on a same
arc of their pencil. This corresponds to $\phi(t)$ having, or not, all
its real roots of the same sign. The simplest way to translate it into
algebraic identities is by using Descartes' law of signs (precisely
Lemma \ref{Descartes:hyperbolic}, because $\phi$ has all its roots
real and non-zero in the considered cases). 

\begin{proposition}
Let $[f],[g]$ be two distinct proper non-empty conics, 
generating a pencil in one of the orbits: {\PI}, {\PIa}, {\PII} {\PIIa}, {\PIII}.
Then $\{[f],[g]\}$ is in the class {\CN} if and only if 
\[
\Phi_{30}\Phi_{12}>0 \wedge \Phi_{03}\Phi_{21}>0
\] 
\end{proposition}

\subsection{Which is inside ?}

Suppose the pair of conics is in one of the classes: {\CIaN},
{\CIIN}, {\CIIaN}, {\CIIIN}, {\CIIIaN}, {\CVN}. Which conic lies inside the
other ? Otherwise stated, for any given class of pairs inside, we want
to characterize the corresponding classes of couples. 

\subsubsection{The antisymmetric invariant solves the problem for pair
  classes {\CIIN}, {\CIIaN}, {\CIIIN}, {\CIIIaN}}

The \emph{antisymmetric invariant} is
\[
\mathcal{A}=\Phi_{30} \Phi_{12}^3-\Phi_{03}\Phi_{21}^3.
\]
First it is homogeneous of even degree, $6$, in $f$, as well as in $g$. 
So its sign depends only on the algebraic conics, not on the quadratic
forms defining them.

Consider again Table \ref{table:orbit_representatives}. Set 
\begin{equation}\label{t:V}
f=f_0 + t_1 g_0, \qquad g=f_0 + t_2 g_0.
\end{equation}
From figures \ref{pencils:partI} and \ref{pencils:partII}, for the
cases {\PIa},{\PII} {\PIIa}, {\PIII}, {\PIIIa}, the inner conic is the one nearer
from $f_0$, that is the one whose parameter ($t_1$ or $t_2$) has
smaller absolute value. For case {\PV}, it is the one with whose
parameter is smaller.
Evaluate the antisymmetric invariant on $(f,g)$.
For {\PII}, {\PIIa}, {\PIII}, {\PIIIa}, we get each time a positive
rational number times 
\[
\left(t_1 t_2 (t_1 - t_2)\right)^2 (t_1^2-t_2^2).
\]
This proves the following proposition. 
\begin{proposition}
Suppose $([f],[g])$ is a couple of distinct proper non-empty conics, such
that $\{[f],[g]\}$  is in class {\CIIN}, {\CIIaN}, {\CIIIN} or {\CIIIaN}.
Then $[f=0]$ lies inside\footnote{only at the neighborhood of
  the double intersection point for class {\CIIN}.} of $[g=0]$ if and
only if $\mathcal{A}(f,g)<0$.
\end{proposition}

For {\CVN}, the evaluation of the antisymmetric invariant gives
zero, and for {\CIaN} it gives the expression 
\[
(t_1^2 -t_2^2)(t_1-t_2)^2 \left(
(t_1+t_2)^2-(t_1 t_2 - 3)^2
\right)
\]
whose sign is not clear.
We need other methods to solve the question in these two cases.

\subsubsection{The antisymmetric covariant solves the problem for class {\CVN}}

Instead of considering the antisymmetric invariant,
we can consider the following \emph{antisymmetric
  covariant}\footnote{This is a quadratic form, and actually the
  antisymmetric invariant of the previous paragraph is its discriminant.}:
\[
\mathcal{B}(f,g)=
\Phi_{12} f - \Phi_{21} g.
\]
We consider its value on $f$, $g$ generating a pencil in orbit {\PV}. It
is enough to look at Levy's representative.
Consider $f$, $g$ as in (\ref{t:V}) for Levy's representative of orbit
{\PV}. Then
\[
\mathcal{B}(xz-y^2+t_1 x^2,xz-y^2+ t_2 x^2)=\frac{t_1-t_2}{4} x^2.
\]
Thus $\mathcal{B}(f,g)$ is a semi-definite quadratic form, negative
when $t_1<t_2$ (that is $[f=0]$ lies inside $[g=0]$) and positive in
the opposite case.
 
For the purpose of calculation, we use that one decides if a
semi-definite quadratic form is negative or positive merely by
considering the sign of the trace of its matrix. Define $T=\operatorname{tr}(\mathcal{B}(f,g))$.
\begin{proposition}
If $f,g$ generate a pencil in orbit {\PV}, then the conic $[f=0]$ lies in inside the conic $[g=0]$ if and only if $T<0$.
\end{proposition}

\subsubsection{Case {\CIaN}}

This case is more difficult than the previous ones.

Suppose $(f,g)$ is in class
{\CIaN}. 
After Figure \ref{pencils:partI}, $\phi(t)$ has three roots of the same sign, two
making $\psi>0$ (conics of the pencil degenerating into isolated points)
and one making $\psi<0$ (conic degenerating into a
double line). Denote them with
$t_1,t_2,t_3$, such that $|t_1|<|t_2|<|t_3|$. Denote also with $\nu$
their common sign (note it is obtained as the sign of $-\Phi_{30}\Phi_{03}$).
The sign of $\Phi_{30}\phi''(t_1)$ is $-\nu$ and the sign of
$\Phi_{30}\phi''(t_3)$ is $\nu$ (because $\Phi_{30}\phi''$ is linear,
with leading coefficient $6\;\Phi_{30}^2$, positive, so it is increasing;
its root lies between $t_1$ and
$t_3$). 

The sign of $\Phi_{30}\phi''(t_2)$ is unknown, denote it with
$\varepsilon$. 

After Figure \ref{pencils:partI}, $[f=0]$ (resp. $[g=0]$) is inside the other \emph{iff} $\psi(t_1)<0$ (resp. $\psi(t_3)<0$).
Thus we have the following table of signs:
\[
\begin{array}{|c|c|c|c|c|}
\cline{3-5}
\multicolumn{2}{c|}{} & t_1  &  t_2  &  t_3  \\
\hline
\multicolumn{2}{|c|}{\Phi_{30}\phi''} & -\nu  & \varepsilon & \nu \\
\hline
\psi   &  [f=0] \text{\ inside}   &  -    &    +        &  +  \\
\cline{2-5}
       &  [g=0] \text{\ inside}   &  +    &    +        &  -  \\
\hline
\hline
\Phi_{30}\phi''\psi   &  [f=0] \text{\ inside}   &  \nu    &    \varepsilon        &  \nu  \\
\cline{2-5}
       &  [g=0] \text{\ inside}   &  -\nu    &    \varepsilon
&  -\nu  \\
\hline
\end{array}
\]
One sees that a Sturm query of $\Phi_{30}\phi''\psi$ for $\phi$ will
give $3$ or $1$ in one case, $-3$ or $-1$ in the other, allowing to
obtain the relative position of the
conics. Precisely, the Sturm queries corresponding to the situations \emph{$[f=0]$
inside \emph{vs.} $[g=0]$ inside} are given by the following table:
\[
\begin{array}{|c|c|c|}
\cline{2-3}
\multicolumn{1}{c|}{} & \nu=+ & \nu=- \\
\hline
\varepsilon=+ & 3\; \text{vs.\ } -1 & -1\; \text{vs.\ } 3 \\
\hline
\varepsilon=- & 1\; \text{vs.\ } -3 & -3\; \text{vs.\ } 1 \\ 
\hline
\varepsilon=0 & 2\; \text{vs.\ } -2 & -2\; \text{vs.\ } 2 \\ 
\hline
\end{array}
\]
Let 
\begin{align*}
Q&=\frac{1}{2}\operatorname{Remainder}(\Phi_{30}\phi''\phi'\psi,\phi)\\
 &= q_2 t^2 + q_1 t + q_0.
\end{align*}
Note that $Q$ can be defined in a simpler way from the already introduced
polynomial $P=\operatorname{Remainder}(\Phi_{30}\phi'\psi,\phi)$, that is:
$Q=\frac{1}{2}\operatorname{Remainder}(\phi''P,\phi)$.

Define $B_i=\sr_i(\phi,Q)$ for $i$ between $0$ and $3$.
Then 
\begin{align*}
B_3&=\Phi_{30}\\
B_2&=q_2\\
B_1&=
\left|
\begin{matrix}
\Phi_{30}  & \Phi_{21} & \Phi_{12} \\
0            &   q_2         &   q_1         \\
q_2          &   q_1         &   q_0  
\end{matrix}
\right|.
\end{align*}
 Finally
\[
B_0=-\Res(\phi,Q).
\]
This last polynomial simplifies.
Using Lemma \ref{res:remainder}, one gets:
\[
\Res(P\phi'',\phi)=- 8 \Phi_{30}\Res(\phi,Q)= 8 \Phi_{30} B_0.
\]
On the other hand, from Lemma \ref{res:product},
\[
\Res(P\phi'',\phi)=\Res(P,\phi)\Res(\phi'',\phi).
\]
From Lemma \ref{res:exchange}, $\Res(P,\phi)=\Res(\phi,P)$, and this
is $-A_0$, which was proved to be equal to:
\[
\Phi_{30}^2 \Res(\psi, \phi) \Disc(\Phi).
\]
Gathering this information, we get that:
\[
B_0=\frac{1}{8} \Phi_{30} \Res(\psi,\phi) \Res(\phi'',\phi) \Disc(\Phi).
\]
It is convenient to remark here that $\Phi_{30}$ divides
$\Res(\phi'',\phi)$.
We will define 
\[
R:=\frac{\Res(\phi'',\phi)}{8 \Phi_{30}}
\]
Thus
\[
B_0= \Phi_{30}^2 \Res(\psi,\phi) R \Disc(\Phi).
\]
From lemmas \ref{res:exchange} and \ref{res:roots}, it comes that 
$\Res(\psi,\phi)=\Res(\phi,\psi)<0$ and $\Res(\phi'',\phi)=-\Res(\phi,\phi'')$ has the sign
$\varepsilon$. Last, $\Disc(\Phi)<0$. Thus $B_0$ has the sign of
$\varepsilon \Phi_{30}$.

The sign sequences $s_1,s_2,s_3,s_4$ giving $3$ or $-3$ are
characterized by $s_1 s_3 >0$ with $s_2 s_4 >0$.

The sign sequences giving $2$ are $+++0$ and $---0$, those giving $-2$
are $+-+0$ and $-+-0$. The first are characterized with respect to the second
by $s_1 s_2>0$.

If $\varepsilon \nu >0$, then $[f=0]$ is inside iff $\Phi_{30} B_1
>0$ with $q_2 \varepsilon \Phi_{30}>0$. If $\varepsilon \nu <0$, then
this characterizes \emph{$[g=0]$ inside}. If $\varepsilon=0$, $[f=0]$ is
inside iff $\nu \Phi_{30} q_2 >0$.

Using that $\varepsilon$ is obtained as the sign of $\Res(\phi'',\phi)$:

\begin{proposition}
Suppose $f,g$ are two non-degenerate quadratic forms 
generating a pencil in orbit {\PIa}. Suppose that their zero locus are
nested, that is $(f,g)$ is in class $N1$. The following are the
necessary and sufficient conditions for $[f=0]$ lies inside $[g=0]$:
\begin{itemize}
\item when $\Phi_{03} R < 0$, it is 
\[
\Phi_{30} B_1 >0 \text{\ and\ } \Phi_{03} q_2  <0;
\]
\item when $\Phi_{03} R >0$, it is 
\[
\Phi_{30} B_1 \leq 0 \text{\ or\ } \Phi_{03} q_2  \leq 0;
\]
\item when $R=0$, it is
\[
\Phi_{03} q_2  <0.
\]
\end{itemize}
\end{proposition}

\subsection{Recapitulation}

Here we display  the explicit definitions of the polynomials appearing
in the description of the rigid isotopy classes. Note that all these
formulas are short: the complicated polynomials express simply in
terms of the less complicated ones. We also display the explicit description of the rigid isotopy classes.

\subsubsection{Formulas}

We will denote the two forms as follows:
\[
\begin{array}{l}
f(x,y,z)=a_{200}x^2+ a_{020} y^2 + a_{002} z^2 + a_{110} xy + a_{101}
xz + a_{011} yz\\
g(x,y,z)=b_{200}x^2+ b_{020} y^2 + b_{002} z^2 + b_{110} xy + b_{101} xz + b_{011} yz.
\end{array}
\]
We will denote similarly the coefficients of
$\tilde{f},\tilde{g},\Omega$ with
$\tilde{a}_{ijk},\tilde{b}_{ijk},\omega_{ijk}$ respectively. One has:
\begin{align*}
\tilde{a}_{200}
=
\left|
\begin{matrix}
a_{020}   &   a_{011}/2 \\
a_{011}/2 &   a_{002} 
\end{matrix}
\right|,
\qquad &
   \tilde{a}_{011} 
   =
   -2\;
   \left|
   \begin{matrix}
   a_{200}   &   a_{110}/2 \\
   a_{101}/2 &   a_{011}/2 
   \end{matrix}
   \right|,
\\
\tilde{a}_{020}
=
\left|
\begin{matrix}
a_{200}   &   a_{101}/2 \\
a_{101}/2 &   a_{002} 
\end{matrix}
\right|,
\qquad &
   \tilde{a}_{101}
   =
   - 2\;\left|
   \begin{matrix}
   a_{020}  & a_{110}/2    \\
   a_{011}/2  & a_{101}/2  
   \end{matrix}
   \right|,
\\
\tilde{a}_{002}
=
\left|
\begin{matrix}
a_{200}   &   a_{110}/2 \\
a_{110}/2 &   a_{020} 
\end{matrix}
\right|,
\qquad &
   \tilde{a}_{110}
   =
   -2\;
   \left|
   \begin{matrix}
   a_{002}   &   a_{011}/2 \\
   a_{101}/2 &   a_{110}/2 
   \end{matrix}
   \right|.
\end{align*}
Similarly the $\tilde{b}_{ijk}$'s are defined from the $b_{ijk}$'s, and 
\begin{align*}
\omega_{200}&=
a_{020}b_{002}+a_{002}b_{020}-a_{011}b_{011}/2,
\\
\omega_{020}&=
a_{002}b_{200}+a_{200}b_{002}-a_{101}b_{101}/2,
\\
\omega_{002}&=
a_{020}b_{200}+a_{200}b_{020}-a_{110}b_{110}/2,
\\
\omega_{011}&=
a_{200}b_{011}+a_{011} b_{200}
- a_{110}b_{101}/2 - a_{101}b_{110}/2,
\\
\omega_{101}&=
a_{020}b_{101}+a_{101} b_{020}
- a_{011}b_{110}/2 - a_{110}b_{011}/2,
\\
\omega_{110}&=
a_{002}b_{110}+a_{110} b_{002}
- a_{011}b_{110}/2 - a_{110}b_{011}/2.
\end{align*}

The (de-homogenized) characteristic form is
\[
\phi(t)=\Phi_{30}t^3+\Phi_{21} t^2 + \Phi_{12} t + \Phi_{03}
       =\Disc(t f+g).
\]
Note that:
\[
\Phi_{30}=
a_{200}\tilde{a}_{200}
+
a_{110}\tilde{a}_{110}
+
a_{101}\tilde{a}_{101},
\]
and
\[
\Phi_{21}=
b_{200}\tilde{a}_{200}
+
b_{002}\tilde{a}_{002}\\
+
b_{020}\tilde{a}_{020}
+
b_{110}\tilde{a}_{110}\\
+
b_{101}\tilde{a}_{101}
+
b_{011}\tilde{a}_{011}.
\]
There are similar formulas for $\Phi_{03}$ and $\Phi_{12}$, by
exchanging $a$ and $b$.

The discriminant of the characteristic form can be obtained as
\[
\Disc(\Phi)=
\frac{1}{81}
\left\vert
\begin{matrix}
3 \Phi_{30}  &  2 \Phi_{21} &  \Phi_{12} & 0\\
 0           &  3 \Phi_{30}  & 2 \Phi_{21}  &  \Phi_{12} \\
\Phi_{21}  & 2 \Phi_{12} & 3 \Phi_{03}  &   0 \\
0            &  \Phi_{21}  & 2\Phi_{12} & 3 \Phi_{03} 
\end{matrix}
\right\vert,
\]
and its Hessian determinant as 
\begin{align}
H
&=
H_{20}t^2 + H_{11} tu + H_{02} u^2
\\
&=
4
\left\vert
\begin{matrix}
3 \Phi_{30} &  \Phi_{21}  \\
 \Phi_{21}  & \Phi_{12}  
\end{matrix}
\right\vert \; t^2
+4
\left\vert
\begin{matrix}
3 \Phi_{30} &  \Phi_{12} \\
 \Phi_{21} & 3 \Phi_{03}  
\end{matrix}
\right\vert\; t\;u
+4
\left\vert
\begin{matrix}
 \Phi_{21} &  \Phi_{12} \\
 \Phi_{12} & 3 \Phi_{03}   
\end{matrix}
\right\vert \; u^2.
\end{align}

The autopolar triangle covariant is:
\begin{multline*}
G=- [\bar{1}\,2\,3] x^3
- [1\,\bar{2}\,3] y^3
- [1\,2\,\bar{3}] z^3
+ ([1\,\bar{1}\,3]+2[\bar{3}\,\bar{2}\,3])x y^2 \\
+ ([1\,2\,\bar{1}]+2[\bar{2}\,2\,\bar{3}])x z^2 
+ ([\bar{2}\,2\,3]+2[\bar{1}\,\bar{3}\,3])y x^2 
+ ([1\,2\,\bar{2}]+2[1\,\bar{1}\,\bar{3}])y z^2\\ 
+ ([\bar{3}\,2\,3]+2[\bar{1}\,2\,\bar{2}])z x^2 
+ ([1\,\bar{3}\,3]+2[1\,\bar{2}\,\bar{1}])z y^2 
+ ([1\,2\,3]+4[\bar{1}\,\bar{2}\,\bar{3}]) xyz.
\end{multline*}
where $[i\,j\,k]$ denote the maximal minors of the matrix 
\[
M=
\left[
\begin{matrix}
\tilde{a}_{200}  & \tilde{a}_{020} & \tilde{a}_{002} & \tilde{a}_{011}
& \tilde{a}_{101} & \tilde{a}_{110} \\
\omega_{200}  & \omega_{020} & \omega_{002} & \omega_{011}
& \omega_{101} & \omega_{110} \\
\tilde{b}_{200}  & \tilde{b}_{020} & \tilde{b}_{002} & \tilde{b}_{011}
& \tilde{b}_{101} & \tilde{b}_{110} \\
\end{matrix}
\right]
\]
whose columns have been labeled $1,2,3,\bar{1},\bar{2},\bar{3}$.

\noindent Denote the coefficients of $\psi$ as follows:
\[
\psi(t)=\Psi_{20}\;t^2+ 2 \; \Psi_{11} \;t + \Psi_{02},
\]
(beware the coefficient of $t$ is $2 \Psi_{11}$)
then
\[
\Psi_{20}=
\tilde{a}_{200}+\tilde{a}_{020}+\tilde{a}_{002},
\]
$\Psi_{02}$ is the corresponding expression with $b$ instead of $a$, and 
\[
\Psi_{11}=
\frac{1}{2} \left( \omega_{200}+\omega_{020}+\omega_{002} \right).
\]
There is also $\mu=\mu_{10}t+\mu_{01}$. Then
\[
\mu_{10}=a_{200}+a_{020}+a_{002}
\]
and $\mu_{01}$ is defined by the corresponding formula with $b$
instead of $a$.

\noindent The polynomial $P$ for the Sturm query of $\psi$ for $\phi$ is
\[
P=\operatorname{Remainder}(\Phi_{30}\phi'\psi,\phi)
=p_2 \;t^2 + p_1\; t + p_0
\]
with 
\begin{align*}
p_2&=
3 \Phi_{30}^2 \Psi_{02} 
- 2 \Phi_{21} \Phi_{30}  \Psi_{11}  
-2 \Phi_{12} \Phi_{30} \Psi_{20} 
+ \Phi_{21}^2 \Psi_{20},\\
p_1&=
2 \Phi_{21} \Phi_{30} \Psi_{02} 
- 4 \Phi_{12} \Phi_{30} \Psi_{11} 
+ \Phi_{12} \Phi_{21} \Psi_{20} 
- 3 \Phi_{03} \Phi_{30} \Psi_{20},\\
p_0&=
\Phi_{12} \Phi_{30} \Psi_{02} 
- 6 \Phi_{03} \Phi_{30} \Psi_{11} 
+ \Phi_{03} \Phi_{21} \Psi_{20}.
\end{align*}
The subresultant $A_1$ is
\[
A_1 =
\left\vert
\begin{matrix}
\Phi_{30}  &  \Phi_{21}  & \Phi_{12}  \\
0          &    p_2      &    p_1     \\
p_2        &    p_1      &    p_0    
\end{matrix}
\right\vert.
\]
The antisymmetric invariant is
\[
\mathcal{A}=\Phi_{30}\Phi_{12}^3-\Phi_{03}\Phi_{21}^3.
\]
and 
the trace of the antisymmetric covariant is
\[
T=
\Phi_{12}\mu_{10}
-
\Phi_{21}\mu_{01}.
\]
The polynomial $Q$ for the Sturm query of $\Phi_{30} \phi'' \psi$ for
$\phi$ is:
\[
Q=\operatorname{Remainder}(P\;\phi'',\phi) 
 =  P\;\phi''-6 p_2 \phi 
 = q_2 \;t^2 + q_1\; t + q_0. 
\]
Its coefficients are
\begin{align*}
q_2 &= 3 p_1 \Phi_{30} - 2 p_2 \Phi_{21}, \\
q_1 &= 3 p_0 \Phi_{30} + p_1 \Phi_{21} -3 p_2 \Phi_{12}, \\
q_0 &= p_0 \Phi_{21} - 3 p_2 \Phi_{03}.
\end{align*}
The subresultant $B_1$ is
\[
B_1=
\left\vert
\begin{matrix}
\Phi_{30}  &  \Phi_{21}  & \Phi_{12}  \\
0          &    q_2      &    q_1     \\
q_2        &    q_1      &    q_0    
\end{matrix}
\right\vert.
\]
The last quantity to consider is 
\[
R=
27 \Phi_{30}^2 \Phi_{03}
+ 2 \Phi_{21}^3 
-6 \Phi_{30} \Phi_{21} \Phi_{12}. 
\]

Each of these expressions is homogeneous in the coefficients of $f$
and as well in the coefficients of $g$. The following table gives
their bi-degree.
\[
\begin{array}{cl@{\qquad}cl@{\qquad}cl}
\tilde{a}_{\alpha} &: (2,0) &
H_{i,j}            &: (2+i,2+j) &
\mathcal{A}        &: (6,6) \\
\omega_{\alpha}    &: (1,1) &  
G                  &: (3,3) &
T                  &: (2,2) \\
\tilde{b}_{\alpha} &: (0,2) & 
\Psi_{i,j}         &: (i,j) &
q_2                &: (8,3) \\
\Phi_{ij}          &: (i,j) &
p_2                &: (6,2) &
B_1                &: (17,8) \\
\Disc(\Phi)        &: (6,6) &
A_1                &: (13,6) &
R                  &: (6,3) 
\end{array}
\]

\subsubsection{The decision procedure}

\paragraph{First step: decide the orbit of pencils.}
Here are the descriptions of the sets of couples of distinct proper
conics generating a pencil in a given orbit.
\begin{align*}
{\PI}:\quad& \Disc(\Phi)<0 \wedge p_2<0 \wedge \Phi_{30} A_1 > 0 \\
{\PIa}:\quad& \Disc(\Phi)<0 \wedge
\left[ p_2 > 0 \vee \Phi_{30}A_1 < 0 \vee \left[ A_1=0 \wedge p_2=0
    \right] \right] \\
{\PIb}:\quad& \Disc(\Phi) > 0  \\
{\PII}:\quad& \Disc(\Phi)=0 \wedge H \neq 0 \wedge G \neq 0 \wedge p_2<0 
\wedge \Phi_{30} A_1 > 0 \\
{\PIIa}:\quad& \Disc(\Phi)=0 \wedge H \neq 0  \wedge G \neq 0 \wedge
\left[ p_2 = 0 \vee \Phi_{30} A_1 < 0 \right] \\
{\PIII}:\quad& \Disc(\Phi)=0 \wedge H \neq 0 \wedge G = 0 \wedge p_2 < 0 \\
{\PIIIa}:\quad& \Disc(\Phi)=0 \wedge H \neq 0 \wedge G = 0 \wedge p_2 > 0 \\
{\PIV}:\quad& H=0 \wedge G\neq 0 \\
{\PV}:\quad& H=0 \wedge G=0.
\end{align*}

\paragraph*{Second step: decide the class of pairs.}
There is only one rigid isotopy class for pairs (class {\CN})
corresponding to each of the orbits of pencils
{\PIb}, {\PIIIa}, {\PIV}, {\PV}.

There are two classes ({\CN} or {\CS}) corresponding to {\PI} {\PIa}, {\PII},
{\PIIa}, {\PIII}. The criterion for being in the class {\CN} is:
\[
\Phi_{30}\Phi_{12}>0 \wedge \Phi_{03}\Phi_{21}>0.
\]

\paragraph{Third step (nested cases): decide which of the conics
  is inside the other.}
The classes of pairs splitting into two classes of couples are:
{\CIaN}, {\CIIaN}, {\CIIIN}, {\CIIIaN}, {\CVN}.

The criteria for $[f=0]$ lies inside $[g=0]$ are the following:
\begin{itemize}
\item {\CIIN}, {\CIIaN}, {\CIIIN}, {\CIIIaN}: $\mathcal{A}<0$.
\item {\CVN}: $T<0$.
\item {\CIaN}:  
\[
\begin{array}{|c|c|c|c|}
\cline{2-4}
\multicolumn{1}{c|}{}&
\multicolumn{3}{|c|}{\text{\ sign\ of\ }\Phi_{03} R}\\
\cline{2-4}
\multicolumn{1}{c|}{}&
-  &  +  &  0 \\
\hline
\text{criterion\ } 
&
\Phi_{30}B_1>0 & \Phi_{03}B_1\leq 0   & \\
\text{for\ }[f=0]& 
\wedge     &  \vee           &  \Phi_{03}q_2<0 \\
\text{inside}&
\Phi_{03}q_2<0 &   \Phi_{03}q_2 \leq 0& \\
\hline
\end{array}
\]
\end{itemize}

\section{Examples and applications}\label{section:examples}

We consider examples and applications for our work. In all of them, we
specialize the above general formulas to pairs of quadratic forms
depending on parameters. We obtain a complicated description of the
partition of the 
parameters space into the subsets corresponding to the isotopy classes.
We then use Christopher Brown's program \emph{SLFQ} of
simplification of  large quantifier-free formulas \cite{SLFQ} to get
simpler descriptions.

\subsection{Two ellipsoids}

We consider two ellipsoids given by the equations (example 2 in \cite{Wang:Wang:Kim}):
\begin{eqnarray*}
x^2+y^2+z^2-25=0,\\
\frac{(x-6)^2}{9}+\frac{y^2}{4}+\frac{z^2}{16}-1=0.
\end{eqnarray*}
We consider then as equations in $x,y$ of two affine conics
depending on a parameter $z$. This corresponds to using a sweeping
plane to explore the two ellipsoids.
We homogenize the equations in $x,y$  with $t$, thus considering:
\begin{eqnarray*}
f=&x^2+y^2+t^2(z^2-25),\\
g=&\frac{(x-6t)^2}{9}+\frac{y^2}{4}+t^2\left(\frac{z^2}{16}-1\right).
\end{eqnarray*}
The quantity $\Disc(\Phi)$ is here $h=49 z^4 + 2516 z^2 - 229376$.
One checks easily that $h$ has two single real roots
$z_0,-z_0$ with $0<z_0<4$.
When the two conics are proper and non-empty, that is when $-4<z<4$,
one finds, using our formulas, that the following classes can occur:
\begin{itemize}
\item[$\bullet$] {\CIaS} when $h>0$, that is $-4<z<-z_0$ or $z_0<z<4$.
\item[$\bullet$] {\CIIaS} when $h=0$, that is $z=\pm z_0$.
\item[$\bullet$] {\CIbN} when $h<0$, that is $-z_0<z<z_0$.
\end{itemize}
The ellipsoids go each through the other.

\subsection{A paraboloid and an ellipsoid}

Our equations, inequations, inequalities  can tell the relative
position of any two conics, not only ellipses, because of the choice
of working in the projective plane.

Thus we can apply also the method of the previous example to any kind
of quadric. In the following example, one considers a paraboloid and an ellipsoid:
\begin{eqnarray*}
4 x^2-4 xy + 2 y^2 - 4 xz +14 x - 6y
+2z^2-10 z +12=0,\\
3 x^2-4 xy+2y^2
-4 xz+16 x + 2 yz - 12 y
+2 z^2-16 z+ 39=0
\end{eqnarray*}
As before, we consider $z$ as a parameter and homogenize the
equations in $x,y$ with $t$, thus considering:
\begin{eqnarray*}
f=&4 x^2-4 xy + 2 y^2 +t(- 4 xz +14 x - 6y)
+t^2(2z^2-10 z +12),\\
g=&3 x^2-4 xy+2y^2
+t(-4 xz+16 x + 2 yz - 12 y)
+t^2(2 z^2-16 z+ 39)
\end{eqnarray*}
We specialize our equations, inequations and inequalities and run
\emph{SLFQ}. Let $z_0=-1/4$ and $z_1<z_2$ be the two roots of $z^2-12 z + 34$.
One finds that $z_0 < z_1$, and $[f=0]$ is proper non-empty when $z>z_0$,
$[g=0]$ is proper non-empty when $z_1<z<z_2$. When both are proper and
non-empty, one finds that the isotopy class is always {\CIaN}. Thus the ellipsoid is inside the paraboloid.

\subsection{Uhlig's canonical forms}

In \cite{Uhlig}, Uhlig presented representatives for the orbits under $GL(n,\R)$ of couples of
quadratic forms 
generating a non-degenerate pencil.

For conics ($n=3$), it follows from Uhlig's presentation that any couple
of conics can be transformed, by means of $PSL(3,\R)$, into one
with associated couple of matrices among:
\[
\begin{array}{l@{,}l@{\quad}l@{\qquad}l@{,}l@{\quad}l}
\left[
\begin{matrix}
1 &  &  \\
  & 1  &  \\
  &   & 
1
\end{matrix}
\right]
&
\left[
\begin{matrix}
\lambda_1 & & \\
&\lambda_2 &    \\
 & & 
\lambda_3
\end{matrix}
\right]
&
{(U_{11})};
&
\left[
\begin{matrix}
1 &  &  \\
  & 1  &  \\
  &   & -1
\end{matrix}
\right]
&
\left[
\begin{matrix}
\lambda_1 & & \\
&\lambda_2 &    \\
 & & -\lambda_3
\end{matrix}
\right]
&
{(U_{12})};
\vspace{2mm}
\\
\left[
\begin{matrix}
  & 1 &  \\
1 &   &  \\
  &   & 
1
\end{matrix}
\right]
&
\left[
\begin{matrix}
  & \lambda_1 & \\
\lambda_1 &  1 & \\
 & & 
\lambda_2
\end{matrix}
\right]
&
{(U_{21})};
&
\left[
\begin{matrix}
  & 1 &  \\
1 &   &  \\
  &   & -1
\end{matrix}
\right]
&
\left[
\begin{matrix}
  & \lambda_1 & \\
\lambda_1 &  1 & \\
 & & -\lambda_2
\end{matrix}
\right]
&
{(U_{22})};
\vspace{2mm}
\\
\left[
\begin{matrix}
  & 1 &  \\
1 &   &  \\
  &   & 
1
\end{matrix}
\right]
&
\left[
\begin{matrix}
 b & a & \\
a &  -b & \\
 & & 
\lambda
\end{matrix}
\right]
&
{(U_{31})};
&
\left[
\begin{matrix}
  & 1 &  \\
1 &   &  \\
  &   & -1
\end{matrix}\right]
&
\left[
\begin{matrix}
 b & a & \\
a &  -b & \\
 & & -\lambda
\end{matrix}
\right]
&
{(U_{32})};
\vspace{2mm}
\\
\multicolumn{6}{c}{
\begin{array}{l@{,}l@{\quad}l}
\left[
\begin{matrix}
  &  &  1\\
  & 1 &  \\
1 &   & 
\end{matrix}
\right]
&
\left[
\begin{matrix}
  &  & \lambda\\
  &  \lambda & 1 \\
\lambda & 1 & 
\end{matrix}\right]
&
{(U_4)}.
\end{array}
}
\end{array}
\]
To which configuration corresponds each of these normal forms ?

We find simple description for the subsets of each parameters space
corresponding to the isotopy classes. 

As an illustration, we show the
result for $U_{21}$. 

\begin{figure}
\centering
      \includegraphics[scale=0.80]{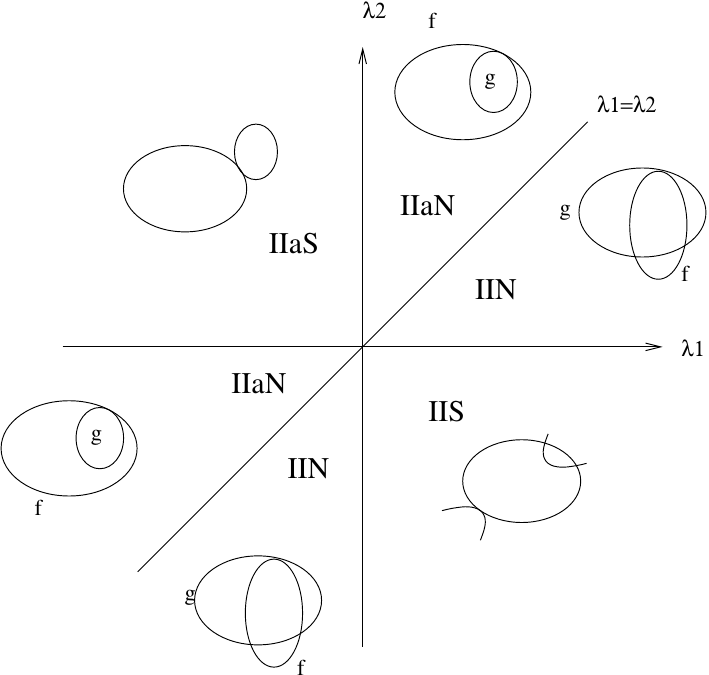}
   \caption{isotopy classes for representatives $U_{21}$.}\label{U21}
\end{figure}

\begin{itemize}
\item $g$ is degenerate when $\lambda_1=0$ or $\lambda_2=0$.
\item the conics are in class {\CVN} when $\lambda_1=\lambda_2 \neq
  0$. In this case, $[g=0]$ lies inside $[f=0]$.
\item in the other cases, the conics are in class {\CIIN}, {\CIIS},
  {\CIIaN} or {\CIIaS}
  as shown in Figure \ref{U21}.
\end{itemize}

\section{Final remarks}

For clarity of the exposition, we have not considered the case
  when one conic, or both conics, are degenerated; but it is easy to
  list the corresponding isotopy classes and describe them with
  equations, inequations and inequalities.

Remark that the polynomials involved in the description of the
classes, specially invariants and covariants, have often very compact
expressions in function of the smaller ones.
Thus they can
be evaluated with substantial saving of arithmetic operations, as was
pointed out in \cite{D:F:M:T}.

The following, more ambitious, step is the classification of couples of
quadrics drawn in $\RP^3$. Hopefully some of the methods developed in
the present paper will be useful in this task, on which we wish to return
in another paper.

It follows from our study that the rigid isotopy classes
  for couples of conics  
  are characterized nearly totally by the behavior of the signature function on the
  pencil generated by the quadratic forms. For the non-generic
  classes, this provides a precise answer to a question formulated in
  \cite{Wang:Krasauskas}. We plan also to develop this point in a
  forthcoming paper with B. Mourrain.

Finally, the reader will find some implementations and complements on
the subject on the author's web page devoted  
to the paper: 
\newline
\url{http://emmanuel.jean.briand.free.fr/publications/twoconics}

\begin{acknowledgements}
The research and the redaction of this paper have been possible thanks
to the successive supports of the European projects GAIA II
(Intersection algorithms for geometry based IT-applications using
approximate algebraic methods); AIM@SHAPE; and of 
 the European RT Network \emph{Real Algebraic and Analytic
   Geometry} (contract No. HPRN-CT-2001-00271). 

The author wants to thank specially Laureano Gonz\'alez-Vega for
introducing him the subject; Bernard Mourrain for fruitful
discussions; Marie-Fran\c{c}oise Coste-Roy and Ioannis Emiris for their
interest; Mercedes Rosas for her careful reading; and the University of Cantabria and INRIA Sophia Antipolis for
their welcoming environment.

The author also wants to thank the anonymous referees for their useful 
comments, and for pointing out the work of Gudkov and Polotovskiy.

Finally, the author is much indebted to the people and institutions
that provide free access to their work on the world wide web: Christopher
Brown for his software SLFQ \cite{SLFQ}, the \emph{project
  Gutenberg}, and the
\emph{Biblioth\`eque Nationale de France} for its project \emph{Gallica}. 
\end{acknowledgements}

\bibliography{deuxconiques}

\end{document}